\documentclass[english]{smfart}
\usepackage{sabbah_Fourier-Hodge}
\begin{document}
\frontmatter
\title[Fourier-Laplace transform]{Fourier-Laplace transform of a variation of~polarized complex Hodge structure}

\author[C.~Sabbah]{Claude Sabbah}
\address{UMR 7640 du CNRS\\
Centre de Math\'ematiques Laurent Schwartz\\
\'Ecole polytechnique\\
F--91128 Palaiseau cedex\\
France}
\email{sabbah@math.polytechnique.fr}
\urladdr{http://www.math.polytechnique.fr/~sabbah}

\thanks{Some part of this work has been achieved as the author was visiting KIAS (Seoul, Korea). The author thanks this institution for providing him with excellent working conditions.}

\begin{abstract}
We show that the Fourier-Laplace transform of a regular holonomic module over the Weyl algebra of one variable, which generically underlies a variation of polarized Hodge structure, underlies itself an integrable variation of polarized twistor structure.
\end{abstract}



\alttitle{Transformation de Fourier-Laplace d'une variation de structure de Hodge complexe polaris\'ee}

\begin{altabstract}
Nous montrons que le transform\'e de Fourier-Laplace d'un module holonome r\'egulier sur l'alg\`ebre de Weyl, sous-jacent g\'en\'eriquement \`a une variation de structure de Hodge polaris\'ee, est lui-m\^eme sous-jacent \`a une variation int\'egrable de structure de twisteur polaris\'ee.
\end{altabstract}

\altkeywords{Fibr\'e plat, variation de structure de Hodge, polarisation, m\'etrique harmonique, $\cD$-module avec structure de twisteur, transformation de Fourier-Laplace}

\maketitle
\mainmatter

\section*{Introduction}

Let $P=\{p_1,\dots,p_r,p_{r+1}=\infty\}$ be a non empty finite set of points on the Riemann sphere $\PP^1$. We will denote by $t$ the coordinate on the affine line $\Afu=\PP^1\moins\{\infty\}$. Let $(V,\nabla)$ be a holomorphic bundle with connection on $\PP^{1\an}\moins P$.\footnote{The exponent $^\an$ means taking the analytic topology.} One can associate to $(V,\nabla)$ a unique holonomic $\Clt$-module $M$ with regular singularities (even at infinity) which is a \emph{minimal extension} on $\Afu$: it is characterized by the fact that its de~Rham complex on $\Afuan$ is $j_*\cV$, if $j:\PP^{1\an}\moins P\hto\Afuan$ denotes the inclusion and $\cV=\ker\nabla$.

The \emph{Laplace transform} $\wh M$ (also called the Fourier-Laplace transform) of the $\Clt$-module $M$ is the $\CC$-vector space $M$ equipped with the following action of the Weyl algebra $\Cltau$: the action of $\tau$ is by $\partial_t$ and that of $\partial_\tau$ is by left multiplication by $-t$ (see \eg \cite{Malgrange91} or \cite[Chap\ptbl V]{Bibi00} for the basic properties of this transformation). We also say that the Laplace transform has kernel $e^{-t\tau}$. In the $\tau$-plane $\Afuhan$, $\wh M$ is a vector bundle with a holomorphic connection $(\wh V,\wh\nabla)$ away from $\tau=0$. It is known that the singularity at $\tau=0$ is regular but the one at infinity is usually irregular (this uses the assumption that $M$ has a regular singularity at $t=\infty$). Recall also that the locally constant sheaf $\wh\cV=\ker\wh\nabla$ can be computed from the locally constant sheaf $\cV=\ker\nabla$ in a cohomological way (see \S\ref{subsec:topolsesqui}), called localized topological Laplace transform.

Let us now assume that $(V,\nabla)$ underlies a variation of polarized complex Hodge structure of some weight~$w\in\ZZ$. We address the following

\begin{question*}
What kind of a structure does the bundle $(\wh V,\wh\nabla)$ associated to the Laplace transform $\wh M$ of $M$ underlie?
\end{question*}

The variation of complex Hodge structure provides $M$ with a good filtration $F_\bbullet M$. In general, there is no way to get from it a good filtration on the Laplace transform. Therefore $(\wh V,\wh\nabla)$ is unlikely to naturally underlie a variation of polarized complex Hodge structure in the classical sense. This is also prevented by the irregular singularity at infinity, according to the regularity theorem of Griffiths and Schmid (\cf \cite{Schmid73}).

One of the main results of this article is Corollary \ref{cor:main}, giving the solution to this question in the following way. We use the language of twistor $\cD$-modules of \cite{Bibi01c}. Let us assume for simplicity that the weight~$w$ is equal to~$0$.

\begin{enumerate}
\item
We extend to $\PP^1$ the main data of the variation, which are only defined on $\PP^1\moins P$. We get $M$ as above, equipped with a good filtration $F_\bbullet M$ and a Hermitian sesquilinear pairing $k:M\otimes_\CC\ov M\to \cS'(\Afu)$ (temperate distributions on the complex plane of the variable $t$).
\item
The \emph{basic correspondence} (\cf Definition \ref{def:basic}) associates to the data~$(M,F_\bbullet M,k)$
\begin{enumerate}
\item
a bundle $G$ with flat connection on the $\tau$-plane away from $\tau=0,\infty$: its analytization was called $(\wh V,\wh\nabla)$ above;
\item
an extension $G_0$ of this bundle across $\tau=\infty$ using the filtration $F_\bbullet M$ by the procedure of \emph{saturation} by $\ptm$;
\item
a sesquilinear pairing between $G$ and $\iota^*G$ (where $\iota$ is the involution $\tau\mto-\tau$) obtained by Fourier transform from $k$.
\end{enumerate}

\item\label{enum:intro3}
The main result (Corollary \ref{cor:main}) says that these data form an \emph{integrable polarized twistor structure of weight~$0$}.

\item\label{enum:intro4}
Moreover, by rescaling the variable $t$ in $\Afu$, we get a corresponding rescaling on the $\tau$-plane, and in this way we get a family of polarized twistor structures of weight~$0$ parametrized by $\CC^*$ (the rescaling factor). We show (\cf Remark \ref{rem:dilat}) that this family is a variation of polarized twistor structure of weight~$0$, with tame behaviour when the rescaling factor (called $1/\tau_o$ in \S\ref{subsec:dilatation}) tends to infinity.

\item
The proof of \eqref{enum:intro3} and \eqref{enum:intro4} is obtained through another interpretation of the objects involved. Indeed, to the original variation of polarized Hodge structure we associate a variation of polarized twistor structure of weight~$0$ in a natural way. Using results of C\ptbl Simpson \cite{Simpson90} and O\ptbl Biquard \cite{Biquard97} as in \cite[Chap\ptbl5]{Bibi01c}, we show that this variation extends as an integrable polarized twistor $\cD$-module of weight~$0$ on $\PP^1$ (we could have used Schmid's classical results \cite{Schmid73}, but we only use the regularity theorem here). We are then in position to apply the main theorem of \cite{Bibi04} (\cf however the erratum to \cite{Bibi04}), saying that the Fourier-Laplace transform of this twistor $\cD$-module is of the same kind. In particular, we get a harmonic metric on $(\wh V,\wh\nabla)$ with a tame behaviour near $\tau=0$.

We then identify the variation of polarized twistor structure that we get on $\Afuh\moins\{0\}$ to the family constructed in \eqref{enum:intro4}. The new point in the proof, taking into account the main theorem in \cite{Bibi04} (with its correction in the erratum to \cite{Bibi04}), is to show that the extension of the original variation of twistor structure is an \emph{integrable} twistor $\cD$-module (in the sense given in \cite[Chap\ptbl7]{Bibi01c}, following the work of C\ptbl Hertling \cite{Hertling01}).

It should be emphasized that, in such a variation, each fibre is naturally equipped with a polarized Hodge structure (of weight~$0$), but the variation does not preserve such a structure, it only preserves the twistor structure, allowing therefore limiting irregular singularities (see \cite[\S7.2]{Bibi01c}).
\end{enumerate}

In \S\ref{sec:Brieskorn}, we use the same ideas to answer a question of C\ptbl Hertling: given a regular function $f:U\to\Afu$ on a smooth affine manifold $U$ of dimension $n+1$, which has only isolated critical points and has a tame behaviour at infinity (\ie is cohomologically tame at infinity, \cf \cite{Bibi96b}), we associate to it the Brieskorn lattice~$G_0$ (a free $\CC[\tau^{-1}]$-module of finite rank); there is also a natural sesquilinear pairing $\wh{C}=\frac{(-1)^{n(n+1)/2}}{(2\pi i)^{n+1}}\wh{\rP}:G_{0|S^1}\otimes_{\cO_{S^1}}\iota^*\ov{G_{0|S^1}}\to\cO_{S^1}$, where the conjugation is taken in the usual sense (the pairing $\wh{\rP}$ will be constructed topologically in \S\ref{sec:Brieskorn}); we then show that $(G_0,G_0,\wh{C})$ corresponds to a polarized integrable twistor structure of weight~$0$. The proof is not obtained by a direct application of the previous results, as the Gauss-Manin system $M=\cH^0f_+\cO_U$ does not usually underlie, generically on~$\Afu$, a variation of polarized Hodge structure (because $f$ is not proper-and-smooth). It underlies a mixed Hodge module in the sense of M\ptbl Saito \cite{MSaito87}. The basic idea is that, under the tameness assumption on $f$, this module differs from a module underlying a variation of polarized Hodge structure only by free $\CC[t]$-modules, which vanish after localized Fourier-Laplace transform, so that the object $(G_0,G_0,\wh{C})$ can also be regarded as associated to a variation of polarized Hodge structure.

When $U$ is a torus $(\CC^*)^{n+1}$ and $f$ is a convenient nondegenerate Laurent polynomial with total Milnor number $\mu$, one defines on the germ $(\CC^\mu,0)$, regarded as the parameter space of a universal unfolding of $f$, a canonical Frobenius structure (\cf \cite{D-S02a}). A consequence of the theorem for $f$ is to endow $(\CC^\mu,0)$ with a positive definite Hermitian metric, satisfying a set of compatibility properties with the Frobenius structure. This is called a $tt^*$-structure in \cite{Hertling01}. Let us notice that, compared to computations made for germs of holomorphic functions in \loccit, the Hermitian form in the case of a Laurent polynomial is positive definite.

\subsubsection*{Acknowledgements}
The results of this article were motivated by many discussions with Claus Hertling, who also carefully read a first version of it. I thank him for all.

The referee pointed out mistakes in \cite{Bibi01c} and in \cite{Bibi04}. While that of \cite{Bibi01c} is easily corrected, that in \cite{Bibi04} needs a detailed correction, which is given in the erratum to \cite{Bibi04}. I gratefully thank the referee for his careful reading of \cite{Bibi01c} and \cite{Bibi04}, as well as for providing me with a key argument in order to complete the proof of Theorem~1 in \cite{Bibi04}.

\section{The basic correspondence}\label{sec:basic}

\Subsection{Fourier-Laplace transform and sesquilinear pairings}\label{subsec:Fsesqui}

\subsubsection*{Conjugation}
Let $X$ be a complex manifold, with structure sheaf $\cO_X$, and let $X_\RR$ denote the underlying $C^\infty$ manifold. We denote by $\ov X$ the manifold $X_\RR$ equipped with the structure sheaf $\cO_{\ov X}\defin \ov{\cO_X}$. The conjugation makes $\cO_{\ov X}$ a $\cO_X$-module and defines a functor, that we call ``conjugation'', from $\cO_X$-modules to $\cO_{\ov X}$-modules.

Given any $\cO_X$-module $\cF$, its conjugate $\cO_{\ov X}\otimes_{\cO_X}\cF$ is denoted by $\ov \cF$: it is a $\cO_{\ov X}$\nobreakdash-module. If $\nabla$ is a flat connection on $\cF$, then $\ov\nabla$ is a flat connection on $\ov\cF$ and $\ker\ov\nabla$ is the local system conjugate to $\ker\nabla$ (corresponding to the conjugate representation of the fundamental group).

Similarly, the notion of conjugation is well-defined for $\cD_X$-modules.

\subsubsection*{Sesquilinear pairings}
A \emph{sesquilinear pairing} on $\cD$-modules $\cM',\cM''$ is a $\cD_X\otimes_\CC\cD_{\ov X}$-linear morphism $\cM'\otimes_\CC\ov{\cM''}\to\Db_{X}$ (the sheaf of distributions on $X_\RR$).

In dimension one, we will also use an affine version of it: let $\Clt$ be the Weyl algebra of the variable $t$ and let $\cS'(\Afu)$ be the Schwartz space of temperate distributions on the complex line. If $M',M''$ are $\Clt$-modules, we will consider sesquilinear pairings $M'\otimes_\CC\ov{M''}\to\cS'(\Afu)$.

If $M'=M''=:M$, we say that a sesquilinear pairing $k:M\otimes_\CC\ov M\to \cS'(\Afu)$ is \emph{Hermitian} if $k(m,\ov n)=\ov{k(n,\ov m)}$ for any $m,n\in M$ (and a similar definition for the sheaf-theoretic analogue).

The sesquilinearity of $k$ allows one to extend $k$ as a sesquilinear pairing from the bicomplex $\DR^\an M^{\prime\an}\otimes_\CC\ov{\DR^\an M^{\prime\prime\an}}$ into the $d',d''$ bicomplex of currents on $\Afu$ (we forget here the behaviour at infinity). As this complex is a resolution of the constant sheaf, one obtains a morphism $k_B$ in the derived category\footnote{The index $_B$ is for ``Betti''.}
\[
k_B:\DR^\an M^{\prime\an}\otimes_\CC\ov{\DR^\an M^{\prime\prime\an}}\to\CC_{\Afuan}.
\]
If we use the notation $\pDR$ (\resp $\pCC$) to denote the de~Rham complex (\resp the constant sheaf) shifted by the dimension of the underlying manifold, this can also be written as
\begin{equation}\label{eq:kB}
k_B:\pDR^\an M^{\prime\an}\otimes_\CC\ov{\pDR^\an M^{\prime\prime\an}}\to\pCC_{\Afuan}[1].
\end{equation}

\subsubsection*{Laplace transform of a $\Clt$-module}
If $M$ is a $\Clt$-module, we denote by $\wh M$ its Laplace transform: this is the $\CC$-vector space $M$ equipped with the following structure of $\Cltau$-module: $\tau$ acts as $\partial_t$ and $\partial_\tau$ acts as $-t$. Given a $\Cltau$-module $N$, we denote by $\iota^*N$ (recall that $\iota$ denotes the involution $\tau\mto-\tau$) the $\CC$-vector space $N$ equipped with the following structure of $\Cltau$-module: $\tau$~acts by $-\tau$ and $\partial_\tau$ acts by $-\partial_\tau$.

The Laplace transform can be obtained in a sheaf theoretic way. We will work on $\PP^1$ and we denote by $\cO_{\PP^1}(*\infty)$ the sheaf of meromorphic functions on $\PP^1$ having pole at most at infinity, so that $\Gamma(\PP^1,\cO_{\PP^1}(*\infty))=\CC[t]$.

We denote by $(*\infty)$ the effect of tensoring (over $\cO_{\PP^1}$) with $\cO_{\PP^1}(*\infty)$ and we call this operation ``localization away from infinity''. In particular $\cD_{\PP^1}(*\infty)$ denotes the sheaf of analytic differential operators localized away from infinity. To a $\Clt$-module $M$ is associated a $\cD_{\PP^1}(*\infty)$-module $\ccM(*\infty)$, that we usually consider as a $\cD_{\PP^1}$-module. Recall that, if $M$ is holonomic, then $\ccM(*\infty)$ is $\cD_{\PP^1}$-holonomic. For such a $M$, let us denote by $p^+\ccM(*\infty)$ its inverse image on $\PP^1\times\Afuh$ (corresponding to $\CC[\tau]\otimes_\CC M$). Let us also denote by $\cE^{-t\tau}$ the free rank-one $\cO_{\PP^1\times\Afuh}(*\infty)$-module with the connection induced by $d-\tau dt-td\tau$. In the following, we will use the notation $p^+\ccM(*\infty)\cE^{-t\tau}$ for the $\cD_{\PP^1\times\Afuh}$-module $p^+\ccM(*\infty)\otimes_{\cO_{\PP^1\times\Afuh}}\cE^{-t\tau}$, that is, the $\cO_{\PP^1\times\Afuh}$-module $p^*\ccM(*\infty)$ with connection twisted by $e^{-t\tau}$. Then
\begin{equation}\label{eq:whM}
\wh M=q_+(p^+\ccM(*\infty)\cE^{-t\tau}),
\end{equation}
 where $q$ is the projection to~$\Afuh$ (see \eg \cite{Malgrange91}).

\subsubsection*{Fourier transform of a sesquilinear pairing}
The Fourier transform $F_t$ with kernel $\exp(\ov{t\tau}-t\tau)\itwopi
dt\wedge d\ov t$ is an isomorphism between $\cS'(\Afu)$ ($t$-plane)
and $\cS'(\Afuh)$ (\hbox{$\tau$-plane}). Given a $2$-form $\psi$ in 
the Schwartz space $\cS(\Afuh)$ (\ie $\psi=\chi(\tau)d\tau\wedge d\ov\tau$ 
with $\chi$ $C^\infty$, rapidly decaying as well as all its 
derivatives when $\tau\to\infty$), we set, for $u\in\cS'(\Afu)$,
\[
\langle F_tu,\psi\rangle\defin\langle u,F_\tau\psi \itwopi
dt\wedge d\ov t\rangle,\quad \text{with }F_\tau\psi=\int 
e^{\ov{t\tau}-t\tau}\psi.
\]
(Recall that $F_\tau\psi$ belongs to $\cS(\Afu)$.) If
$k:M'\otimes_\CC\ov{M''}\to\cS'(\Afu)$ is a sesquilinear pairing, we
denote by $F_tk$ the composition $F_t\circ k$ of $k$ with the Fourier transform of
temperate distributions. Then $F_tk$ becomes a sesquilinear pairing
\[
F_tk:\wh{M'}\otimes_\CC\ov{\iota^+\wh{M''}}\to \cS'(\Afuh).
\]
(the $\iota^+$ is needed as we use the kernel $e^{\ov{t\tau}}$ for the
Laplace transform of $\ov M''$, not $e^{-\ov{t\tau}}$). Let us notice that,
at this stage, $k$ can be recovered from $F_tk$ by composing with the
inverse Fourier transform.

\subsubsection*{The case of holonomic $\Clt$-modules with regular singularity at infinity}
Let us now assume that $M',M''$ are holonomic $\Clt$-modules which have a regular singularity at infinity. Then $\wh{M'},\wh{M''}$ have singularities at $\tau=0$ and $\tau=\infty$ only. Denote by $\wh V',\wh V''$ the holomorphic vector bundles with connection $\wh\nabla$ obtained by restricting $\wh{M'},\wh{M''}$ to $\tau\neq0$, and by $\wh\cV',\wh\cV''$ the corresponding local systems $\ker\wh\nabla$.

The sesquilinear pairing $F_tk$ induces a sesquilinear pairing
\begin{equation}\label{eq:Ftk}
F_tk:\wh V'\otimes_\CC \ov{\iota^+\wh V''}\to\cC^\infty_{\{\tau\neq0\}},
\end{equation}
the datum of which is equivalent to that of a sesquilinear pairing of local systems
\begin{equation}\label{eq:FtkB}
(F_tk)_B:\wh\cV'\otimes_\CC\ov{\iota^{-1}\wh\cV''}\to \CC_{\{\tau\neq0\}},
\end{equation}
and, denoting by $S^1$ the circle $\module{\tau}=1$, it is equivalent to the datum of a sesquilinear pairing
\begin{equation}\label{eq:sesqui}
(F_tk)_B:\wh\cV'_{|S^1}\otimes_\CC\ov{\iota^{-1}\wh\cV''_{|S^1}}\to \CC_{S^1}.
\end{equation}

\skippointrait
\begin{exemples}\label{ex:basics}
\begin{enumerate}
\item\label{ex:basics1}
Let us first assume that $M'=M''=M$ is equal to $\Clt/(t-c)$ for some $c\in\CC$. Denote by $m$ the class of $1$ in $M$ and by $\delta_c$ the distribution which satisfies, for any $C^\infty$ function $\varphi$, the equality $\langle \delta_c,\varphi \itwopi\,dt\wedge d\ov t\rangle=\varphi(c)$. Then, up to a constant, we have $k(m,\ov m)=\delta_c$ and $F_tk(m,\ov m)=e^{\ov{c\tau}-c\tau}$.
\item\label{ex:basics2}
Let us now assume that $M'=M''=M$ is equal to $\Clt/(t\partial_t-\alpha)$ for some $\alpha\in{}]-1,0[$. Denote by $m$ the class of $1$ in $M$. Then, up to a constant, $k(m,\ov m)= \module{t}^{2\alpha}$ and
$F_tk(m,\ov m)= 
\frac{\Gamma(\alpha+1)}{\Gamma(-\alpha)}\module{\tau}^{-2(\alpha+1)}$.
\end{enumerate}
\end{exemples}

\subsubsection*{Computation using direct images}
In order to compare with the topological Fourier transform (see below), it will be convenient to have another formulation of the Fourier transform. Recall that $\wh M$ can be regarded as $q_+(p^+M\otimes\cE^{-t\tau})$ (see the diagram below for the notation $p,q$), that is, the cokernel of the injective morphism
\[
\CC[\tau]\otimes_\CC M\To{\nabla_t-\tau dt}\CC[\tau]\otimes_\CC M\otimes dt
\]
[where $\nabla_t$ is the connection relative to $t$ only] \emph{via} the map $\sum \tau^im_i\, dt\mto\sum_i(\partial_t)^im_i\in M$, and the action of $\partial_\tau$ is induced by that of $\partial_\tau-t$ on $\CC[\tau]\otimes_\CC M$. We consider the pairing
\begin{multline}\label{eq:Fouriertau}
\textstyle(\sum_i\tau^im'_i)dt\otimes\ov{(\sum_j\tau^jm''_j)dt}\\
\longmapsto
\Big[\cS(\Afuh)^{(1,1)}\ni\psi\mapsto\sum_{i,j}\langle 
k(m'_i,\ov{m''_j}),F_\tau(\tau^i\ov\tau^j\psi)dt\wedge d\ov t\rangle\Big].
\end{multline}
Let us notice that it vanishes if one of the terms belongs to $\im (\partial_t-\tau)$, hence naturally defines a sesquilinear pairing between the cokernels of $\partial_t-\tau$ with values in $\cS'(\Afuh)$, that we denote by $\wh k$. The following is then clear:

\begin{lemme}\label{lem:itwopi}
We have $F_tk=\itwopi\wh k$.\qed
\end{lemme}

\Subsection{Topological Fourier-Laplace transform and sesquilinear pairings}\label{subsec:topolsesqui}
With the supplementary assumption that $M',M''$ have regular singularity everywhere, the sesquilinear pairing $\wh k_B=-2\pi i(F_tk)_B$ can be obtained from $k_B$ (defined by \eqref{eq:kB}) in a topological way. We will explain here the relationship between these pairings.

\subsubsection*{Definition of the topological Laplace transform of sheaves}
In this paragraph, we use the analytic topology of $\Afu$, $\PP^1$ or $\Afuh$, so we do not indicate it by the exponent `$\an$'. Although the following construction holds over $\Afuh$, we only use it (and therefore explain it) out of $\tau=0$, that is, over $\Afuhs$.

We denote by $\ell$ the open inclusion $\Afu\hto\PP^1$ or $\Afu\times\Afuhs\hto\PP^1\times\Afuhs$ (and, as above, we have $j:\PP^1\moins P\hto\Afu$). Let $\epsilon:\wt\PP^1\to\PP^1$ be the real oriented blowing-up of $\infty\in\PP^1$ ($\wt\PP^1$ is topologically a disc). We also denote by $\epsilon$ the induced map $\wt\PP^1\times\Afuhs\to\PP^1\times\Afuhs$ and we denote by $L^{\prime+}\subset\epsilon^{-1}(\infty)\simeq S^1\times\Afuhs$ the closed subset $\reel(\tau e^{i\arg t})\geq0$ and by $L^{\prime-}$ its complement in $\wt\PP^1\times\Afuhs$. We will consider the commutative diagram
\[
\xymatrix@C=1.5cm{
\Afu\times\Afuhs \ar@<-1mm>`u[r]`[rr]^-{\wt\ell}[rr]
\ar@{_{ (}->}[drr]_{\ell}\ar@<-.5mm>@{^{ (}->}[r]^-{\alpha'}&
L^{\prime-}\ar@<-.5mm>@{^{ (}->}[r]^-{\beta'}&\wt\PP^1\times\Afuhs\ar@<-1mm>[d]^-{\epsilon}\ar[ddr]^(.3){\wt q}
\ar[ddl]_(.3){\wt p}|(.4)\hole\\
&&\PP^1\times\Afuhs\ar[dl]_{p}\ar[dr]^{q}\\
&\PP^1&&\Afuhs
}
\]

Let $\cG$ be a complex of sheaves on $\Afu$. Recall (see \eg \cite{Malgrange91}) that the topological Laplace transform of $\cG$ with kernel $e^{-t\tau}$ (restricted to $\Afuhs$) is the complex
\[
\wh\cG\defin\bR\wt q_*\Big[\beta'_!\bR\alpha'_*p^{-1}\cG\Big][1],
\]
where we still denote by $p$ the projection $\Afu\times\Afuhs\to\Afu$. This definition can be simplified if we assume that $\cG$ is a constructible sheaf: then $\cG$ is a local system near $\infty$, and we have
$\bR\alpha'_*p^{-1}\cG=\alpha'_*p^{-1}\cG$. Moreover,
$\beta'_!\alpha'_*p^{-1}\cG$ commutes with the restriction to
$\tau_o\in\Afuhs$, that is, denoting by $L^{\prime-}_{\tau_o}$ the
intersection $L^{\prime-}\cap\wt\PP^1\times\{\tau_o\}$ and by
$\alpha'_{\tau_o},\beta'_{\tau_o}$ the corresponding inclusions, we
have
$(\beta'_!\alpha'_*p^{-1}\cG)_{|\wt\PP^1\times\{\tau_o\}}=\beta'_{\tau_o,!}
\alpha'_{\tau_o,*}\cG$. By base change for a proper morphism, we then
have
\[
\cH^j\wh\cG_{\tau_o}=\bR^{j+1}\wt q_*\Big[\beta'_!\bR\alpha'_*p^{-1}\cG\Big]_{\tau_o}=
H^{j+1}\big(\wt\PP^1,\beta'_{\tau_o!}\alpha'_{\tau_o*}\cG\big).
\]
If $\Phi_{\tau_o}$ denotes the family of closed sets of $\Afu$, the closure of which in $\wt\PP^1$ does not cut $L_{\tau_o}^{\prime+}$, we have by definition
\[
\bH^{j+1}\big(\wt\PP^1,\beta'_{\tau_o!}\alpha'_{\tau_o*}\cG\big)=
H^{j+1}_{\Phi_{\tau_o}}\big(\Afu,\cG\big).
\]

Let us now assume that $\cG$ is a $\CC$-perverse sheaf\footnote{We refer for instance to \cite{Dimca04} for basic results on perverse sheaves; recall that the constant sheaf supported at one point is perverse, and a local system shifted by one is perverse, see \eg \cite[Ex\ptbl5.2.23]{Dimca04}. Let us notice that, in this paragraph and in the next one, one can work with $\QQ$-perverse sheaves.}. Then $\bH^{j+1}\big(\wt\PP^1,\beta'_{\tau_o!}\alpha'_{\tau_o*}\cG\big)=0$ for $j\neq-1$\footnote{This can be proved as follows: using the structure theorem for perverse sheaves on $\Afu$, on reduces to the case of a sheaf supported on some $p_j$ (trivial), and to the case of $j_*\cV[1]$, where $\cV$ is a local system on $\Afu\moins\{p_1,\dots,p_r\}$; clearly, there is no $H^0_{\Phi_{\tau_o}}(\Afu,j_*\cV)$ and, by duality, there is no $H^2_{\Phi_{\tau_o}}(\Afu,j_*\cV)$.}. In other words, the complex $\wh\cG$ has cohomology in degree $-1$ only. Up to a shift by $-1$, it is a local system on $\Afuhs$ with fiber at $\tau_o$ equal to $\bH^0_{\Phi_{\tau_o}}\big(\Afu,\cG\big)$. Hence $\wh\cG$ is a smooth perverse sheaf on $\Afuhs$.

The Laplace transform with kernel $e^{t\tau}$ is defined similarly, using $L^{\prime\prime\pm}$ obtained with $\reel(\tau e^{i\arg t})\leq0$, by replacing $\alpha',\beta'$ with $\alpha'',\beta''$.

Let us notice that $L^{\prime-}\cap L^{\prime\prime-}=\Afu\times\Afuhs$.

\subsubsection*{Topological Fourier transform of a sesquilinear pairing}
Let $\cG',\cG''$ be $\CC$-perverse sheaves on $\Afu$. Recall that we denote by $\pCC$ the constant sheaf shifted by the dimension of the underlying manifold.

Let us assume that we are given a morphism (in the derived category of bounded complexes with constructible cohomology)
$k_B:\cG'\otimes_\CC\ov{\cG''}\to\pCC_{\PP^1}[1]$. We then get a
morphism
\begin{equation}\label{eq:sheafsesqui}
\beta'_!\alpha'_*p^{-1}\cG'\otimes_\CC\ov{\beta''_!\alpha''_*p^{-1}\cG''}\to \wt\ell_!\pCC_{\Afu\times\Afuhs}.
\end{equation}
For any $\tau_o\in\Afuhs$, we obtain a sesquilinear pairing
\begin{equation}\label{eq:kBo}
\wh{k_B}_{\tau_o}:\bH^0_{\Phi_{\tau_o}}\big(\Afu,\cG'\big)\otimes_\CC 
\bH^0_{\Phi_{-\tau_o}}\big(\Afu,\ov{\cG''}\big)\to H^2_c(\Afu,\CC)\simeq\CC
\end{equation}
using that a closed set in $\Afu$ is both in $\Phi_{\tau_o}$ and $\Phi_{-\tau_o}$ iff it is compact. This pairing defines a sesquilinear pairing at the sheaf level between local systems:
\begin{equation}\label{eq:sheafFLsesqui}
\wh{k_B}:\wh{\cG'}_{|S^1}[-1]\otimes_\CC\ov{\iota^{-1} 
\wh{\cG''}_{|S^1}}[-1]\to\CC_{S^1}.
\end{equation}

\subsubsection*{Computation of $\wh\cF$}
We keep notation as above and we set $M=M',M''$ and $\cF=\DR^\an M[1]$, which is a $\CC$-perverse sheaf, if $\DR$ denotes the usual de~Rham functor. Let us assume that we are given a sesquilinear pairing $k$. In order to compare $\wh{k_B}$ and $\wh k_B$, we need to consider a space where both pairings are defined simultaneously, and to sheafify the construction of $\wh k$ on this space. We will work on $\wt\PP^1\times\Afuhs$.

Let us first recall the natural identification of local systems on $\Afuhs$
\begin{equation}\label{eq:identtauo}
\wh\cF[-1]\simeq\wh\cV,
\end{equation}
when $M$ has a regular singularity at infinity. We have
\[
\wh\cV=\DR^\an(\wh M_{|\tau\neq0})=\DR^\an(\wh V,\wh\nabla).
\]

By restricting \eqref{eq:whM} to $\tau\neq0$, we find $\wh M_{|\tau\neq0}=q_+(p^+\ccM(*\infty)\cE^{-t\tau})$, and we have, as $q$ is proper (so that we can use GAGA relative to $\PP^1$),
\begin{equation}\label{eq:pDR}
\wh\cV=\bR q_*\big(\DR^\an(p^+\ccM(*\infty)\cE^{-t\tau})\big)[1].
\end{equation}

If $A$ is a subset of $\wt\PP^1\times\Afuhs$, we denote by $A^*$ its intersection with $\Afu\times\Afuhs$. Let $K$ be a compact set in $\wt\PP^1\times\Afuhs$ and let $\varphi$ be a $C^\infty$ function on $\Afu\times\Afuhs$ supported in $K^*$. One sets
\[
\cN_p(\varphi)=\sum_{\substack{k\leq p\\ \module{\alpha}\leq p}}\big\Vert\module{t}^{2k}\partial^\alpha\varphi\big\Vert_{L^\infty},
\]
where $\alpha$ is a multi-index indicating derivations with respect to $t,\ov t,\tau,\ov\tau$.

In order to give a realization of the complex \eqref{eq:pDR}, we introduce the following sheaves on $\wt\PP^1\times\Afuhs$. We will not distinguish between distributions and currents, by fixing the volume forms $\itwopi\,dt\wedge d\ov t$ and $\itwopi\,d\tau\wedge d\ov\tau$.

\begin{itemize}
\item
The sheaf $\cE^{<0}_{\wt\PP^1\times\Afuhs,c}$ ($C^\infty$ functions on $\Afu\times\Afuhs$ with compact support in $\wt\PP^1\times\Afuhs$, with rapid decay, as well as all their derivatives, along $\epsilon^{-1}(\infty)$): for any open set $U\subset \wt\PP^1\times\Afuhs$ and any compact set $K\subset U$, $\varphi \in\cE^{<0}_{\wt\PP^1\times\Afuhs,K}(U)$ iff $\varphi\in\cC^\infty(U^*)$ has support in $K^*$ and, for any $p\in\NN$, $\cN_p(\varphi)<+\infty$.
\item
The sheaf $\cEmod_{\wt\PP^1\times\Afuhs}$ ($C^\infty$ functions on $\Afu\times\Afuhs$ having moderate growth, as well as all their derivatives, along $\epsilon^{-1}(\infty)$): a section $\varphi\in\cEmod_{\wt\PP^1\times\Afuhs}(U)$ is a $C^\infty$ function on $U^*$ such that, for any compact set $K\subset U$ and any $p\in\NN$, there exists an integer $N=N(K,p)$ such that $\big\Vert\module{t}^{-2N}\varphi\big\Vert_{C^p,K}<+\infty$.
\item
The sheaf $\cEmod_{\wt\PP^1\times\Afuhs,c}$: same as above, with compact support. So $\varphi$ has support in some $K$ and for any $p$ there exists $N=N(p)$ such that $\big\Vert\module{t}^{-2N}\varphi\big\Vert_{C^p}<+\infty$.
\item
The sheaf $\cEmodan_{\wt\PP^1\times\Afuhs}$: the subsheaf of $\cEmod_{\wt\PP^1\times\Afuhs}$ of functions which are holomorphic with respect to the $\tau$ variable (\ie killed by $\partial_{\ov\tau}$). Using a Cauchy-type argument with respect to $\tau$, it is enough, to control the moderate growth of the derivatives with respect to $t,\ov t$.
\item
The sheaf $\cAmod_{\wt\PP^1\times\Afuhs}$: the subsheaf of $\cEmod_{\wt\PP^1\times\Afuhs}$ of functions which are holomorphic on $\Afu\times\Afuhs$. Using a Cauchy-type argument, it is enough to control the moderate growth of the function itself, not its derivatives.
\item
The sheaf $\Dbmod_{\wt\PP^1\times\Afuhs}$: for any open set $U$ of $\wt\PP^1\times\Afuhs$, the space $\Dbmod(U)$ is the space of linear forms on $\cE^{<0}_{\wt\PP^1\times\Afuhs,c}(U)$ such that, for any compact set $K\subset U$, there exists $p\in\NN$ and $C\geq0$ such that, for any $\varphi\in\cE^{<0}_{\wt\PP^1\times\Afuhs,K}(U)$, one has $\module{\langle u,\varphi\rangle}\leq C\cN_p(\varphi)$.
\item
The sheaf $\Db^{<0}_{\wt\PP^1\times\Afuhs}$: for any open set $U$ of $\wt\PP^1\times\Afuhs$, the space $\Db^{<0}(U)$ is the space of linear forms on $\cEmod_{\wt\PP^1\times\Afuhs,c}(U)$ such that, for any compact set $K\subset U$ and for any integer $N\geq0$, there exists an integer $p=p(K,N)\geq 0$ and a number $C=C(K,N)>0$ such that, for any $\varphi\in\cEmod_K(U)$, one has $\big\vert\langle u,\module{t}^{2N} \varphi\rangle\big\vert\leq C\norme{\varphi}_{C^p}$. Let us notice that $\Db^{<0}(U)\subset\Db(U)$.
\end{itemize}

\begin{proposition}\label{prop:DolbeaultDR}
The previous sheaves are stable by the derivations $\partial_t,\partial_{\ov t},\partial_\tau,\partial_{\ov\tau}$. Moreover,
\begin{enumerate}
\item\label{prop:DolbeaultDR1}
We have $\epsilon_*\cAmod_{\wt\PP^1\times\Afuhs}=\cO_{\PP^1\times\Afuhs}(*\infty)$ and $\bR^j\epsilon_*\cAmod_{\wt\PP^1\times\Afuhs}=0$ for $j\geq1$;
\item\label{prop:DolbeaultDR2}
The Dolbeault complex $\cE^{\rmod,\an,(0,\cbbullet)}_{\wt\PP^1\times\Afuhs}$ is a resolution of $\cAmod_{\wt\PP^1\times\Afuhs}$, that is, $\partial_{\ov t}:\cEmodan_{\wt\PP^1\times\Afuhs}\to\cEmodan_{\wt\PP^1\times\Afuhs}$ is onto and its kernel is $\cAmod_{\wt\PP^1\times\Afuhs}$.
\item\label{prop:DolbeaultDR3}
The complexes $\DR\Dbmod_{\wt\PP^1\times\Afuhs}$ and $\DR\Db^{<0}_{\wt\PP^1\times\Afuhs}$ are resolutions of $\CC_{\wt\PP^1\times\Afuhs}$ and $\wt\ell_!\CC_{\Afu\times\Afuhs}$ respectively.
\end{enumerate}
\end{proposition}

\noindent
[Let us notice that the moderate Dolbeault complex can be computed with $\partial_{\ov t}$ because, if $t'$ is a local coordinate on $\PP^1$ at $\infty$, we have $\partial_{\ov{t'}}=-\ov t^2\partial_{\ov t}$, and multiplication by $\ov t$ is an isomorphism on $\cEmodan_{\wt\PP^1\times\Afuhs}$.]

\begin{proof}[Indication of proof]
For the first point, see \eg \cite[Cor\ptbl II.1.1.8]{Bibi97}. The second point is analogous to Prop\ptbl II.1.1.7 in \loccit, as well as the third point.
\end{proof}

We can compute \eqref{eq:pDR} at the level of $\wt\PP^1$: by the projection formula and Lemma \ref{prop:DolbeaultDR}\eqref{prop:DolbeaultDR1}, the right-hand term of \eqref{eq:pDR} is isomorphic to
\[
\bR \wt q_*\big(\DR^{\rmod}(p^+\ccM(*\infty)\cE^{-t\tau})\big)[1],
\]
where $\DR^{\rmod}$ is the de~Rham complex of
\[
\cAmod_{\wt\PP^1\times\Afuhs}\hspace*{-1mm}\ootimes_{\epsilon^{-1}\cO_{\wt\PP^1\times\Afuhs}}\hspace*{-5mm}\epsilon^{-1}\big[p^+\ccM(*\infty)\cE^{-t\tau}\big].
\]

\begin{lemme}\label{lem:compFderham}
The complex $\DR^{\rmod}(p^+\ccM(*\infty)\cE^{-t\tau})[1]$ is a resolution of $\beta'_!\alpha'_*p^{-1}\cF$.
\end{lemme}

\begin{proof}
Analogous to \cite[Appendix~A]{Malgrange91}.
\end{proof}

The identification \eqref{eq:identtauo} is obtained by using Lemma \ref{lem:compFderham}:
\begin{equation}\label{eq:Rqstar}
\bR \wt q_*\big(\DR^{\rmod}(p^+\ccM(*\infty)\cE^{-t\tau})\big)[1]=\bR \wt 
q_*(\beta'_!\alpha'_*p^{-1}\cF)=\wh\cF[-1].
\end{equation}

Using Dolbeault Lemma \ref{prop:DolbeaultDR}\eqref{prop:DolbeaultDR2}, the $\wt q_*$-acyclicity of $\cEmodan_{\PP^1\times\Afuhs}$ and the $\cO_{\PP^1\times\Afuhs}$-flatness of $p^+\cM(*\infty)\cE^{-t\tau}$ near $\{\infty\}\times\Afuhs$, we find that the previous complex~is
\begin{equation}\label{eq:qstar}
\wt q_*\big(\cE^{\rmod,\an,1+\cbbullet}_{\wt\PP^1\times\Afuhs}\hspace*{-3mm} \ootimes_{\epsilon^{-1}\cO_{\wt\PP^1\times\Afuhs}}\hspace*{-3mm} \epsilon^{-1}[p^+\ccM(*\infty)\cE^{-t\tau}]\big).
\end{equation}

\Subsubsection*{Comparison with the analytic Fourier transform}

\begin{proposition}\label{prop:kBkB}
Under the identification \eqref{eq:identtauo} (and its complex
conjugate for $\iota^+\wh M''$), we have $\wh{k_B}=\wh k_{B}$.
\end{proposition}

\begin{proof}
We will sheafify below the construction of $\wh k$. Let us begin with
a basic fact. Let $\eta=\chi(t,\tau)d\tau\wedge d\ov\tau$ be a 
$2$-form, with $\chi\in\cEmod_c(\wt\PP^1\times\Afuhs)$ having
support in $K$. We denote by $F_\tau(\eta)=\int_{\Afuh} 
e^{\ov{t\tau}-t\tau}\eta$ its Fourier transform
relative to $\tau$. It is a function of $t$.

\begin{lemme}
With these assumptions, for all integers $p\geq0$ and $N\geq0$, there exist $q=q(K,p,N)\in\NN$ and $C=C(K,p,N)>0$ such that
\[
\cN_{p}\big(F_\tau(\eta)\big)\leq C\big\Vert\module{t}^{-2N} \chi\big\Vert_{C^{p+q}}.
\]
\end{lemme}

\begin{proof}
This is a variant of the fact that 
$t\mto\int_{\Afuh}e^{\ov{t\tau}-t\tau}\psi(\tau)$ is in the Schwartz class 
when $\psi$ has compact support.
\end{proof}

\begin{lemme}\label{lem:Fourierlocal}
Let $U$ be an open set of $\wt\PP^1$, $\wh U$ be an open set of
$\Afuhs$ and $u$ be a moderate distribution on $U$ (relative to
$\epsilon^{-1}(\infty)\cap U$). Then the correspondence
$\cEmod_c(U\times\wh U)d\tau\wedge d\ov\tau\ni\eta\mto\langle 
u,F_\tau(\eta)\rangle$
defines an element of $\Db^{<0}(U\times\wh U)$.
\end{lemme}

\begin{proof}
As $u$ is moderate, for any compact set $L$ of $U$ there exist $p$ and
a constant $C>0$ such that, for any function $\varphi\in\cE^{<0}_L(U)$,
one has $\module{\langle u,\varphi \itwopi dt\wedge d\ov t\rangle}\leq 
C\cN_p(\varphi)$. One can
take $\varphi=F_\tau(\eta)$ and $L=\text{projection of }K$; $\varphi$ has
rapid decay, as well as all its derivatives, after the previous lemma,
and one gets, for any $N\geq0$,
\[
\module{\langle u,F_\tau(\eta)\rangle}\leq 
C'(K,p,N)\big\Vert\module{t}^{-2N}\chi\big\Vert_{C^{p+2N}}.\qedhere
\]
\end{proof}

The sesquilinear pairing $k:M'\otimes\ov{M''}\to\cS'(\Afu)$ can be sheafified to give a pairing $k:\ccM'(*\infty)\otimes\ov{\ccM''(*\infty)}\to \Dbmod_{\PP^1}$. For any $a,b\in\{0,1,2\}$, we define a pairing $\wt{k_B}$ between sections of 
$\cAmod_{\wt\PP^1\times\Afuhs}\otimes\epsilon^{-1}\big[\Omega_{\PP^1\times\Afuhs}^a\otimes p^+\ccM'(*\infty)\cE^{-t\tau}\big]$ and sections of $\ov{\cAmod_{\wt\PP^1\times\Afuhs}\otimes\epsilon^{-1}\big[\Omega_{\PP^1\times\Afuhs}^b\otimes p^+\ccM''(*\infty)\cE^{t\tau}\big]}$ with values in the sheaf of $a+b$-currents $\Db_{\wt\PP^1\times\Afuhs}^{<0,a+b}$ by setting, for any $4-(a+b)$-form $\eta$ with coefficients in $\cE^{\rmod}_{\wt\PP^1\times\Afuhs,c}$ with support in the open set where the sections are defined,
\[
\Big\langle \wt{k_B}\big({\textstyle\sum_i\psi'_i\otimes m'_i,\sum_j\ov{\psi''_j\otimes m''_j}}\big), \eta\Big\rangle=\sum_{i,j}\Big\langle k(m'_i,\ov{m''_j}),\int_{\Afuhs}e^{\ov{t\tau}-t\tau} \psi'_i\wedge\ov{\psi''_j}\wedge\eta\Big\rangle.
\]
We note that the right-hand term is meaningful because of Lemma \ref{lem:Fourierlocal}.

\begin{lemme}
The pairing $\wt{k_B}$ induces a pairing of bicomplexes
\[
\DR^{\rmod}(p^+\ccM'(*\infty)\cE^{-t\tau})\ootimes_\CC\ov{\DR^{\rmod}(p^+\ccM''(*\infty)\cE^{t\tau})}\to \Db_{\wt\PP^1\times\Afuhs}^{<0,(\cbbullet,\cbbullet)}
\]
and induces \eqref{eq:sheafsesqui}, as obtained from \eqref{eq:kB}, at the level of the associated simple complexes.
\end{lemme}

\begin{proof}[Sketch of proof]
The first point follows from the sesquilinearity of $k$. By Lemma
\ref{lem:compFderham}, the left-hand term is a resolution of
$\beta'_!\alpha'_*p^{-1}\cF'[-1]\otimes_\CC\beta''_!\alpha''_*p^{-1}\ov{\cF''}[-1]$,
and by Proposition \ref{prop:DolbeaultDR}\eqref{prop:DolbeaultDR3},
the right-hand term is a resolution of
$\wt\ell_!\CC_{\Afu\times\Afuhs}$. As, in any case, the morphism
induced by $\wt{k_B}$ coincides with \eqref{eq:sheafsesqui} along
$\epsilon^{-1}(\infty)$ (where both are zero), it is enough to show
the coincidence locally on $\Afu\times\Afuhs$, where the result is standard. \end{proof}

In order to compute analytically the pairing $\wh{k_B}$, we resolve the de~Rham complexes above with coefficients in $\cA^{\rmod}_{\wt\PP^1\times\Afuhs}$ by $C^\infty$ de~Rham complexes with coefficients in $\cE^{\rmod,\an}_{\wt\PP^1\times\Afuhs}$ in order to get $\wt q$-acyclicity. Then, $\wh{k_B}$ is obtained by applying $\wt q_*$ to the pairing
\begin{multline}\label{eq:wtwhk}
\Big(\cE^{\rmod,\an,1}_{\wt\PP^1\times\Afuhs}\hspace*{-2mm} \ootimes_{\epsilon^{-1}\cO_{\wt\PP^1\times\Afuhs}}\hspace*{-2mm} \epsilon^{-1} \big[p^+\ccM'(*\infty)\cE^{-t\tau}\big]\Big)\\
\otimes
\ov{\Big(\cE^{\rmod,\an,1}_{\wt\PP^1\times\Afuhs}\hspace*{-2mm} \ootimes_{\epsilon^{-1}\cO_{\wt\PP^1\times\Afuhs}}\hspace*{-2mm}\epsilon^{-1}\big[p^+\ccM''(*\infty)\cE^{t\tau}\big]\Big)}\\
\smash{\To{\wt{k_B}}}\,\Db^{<0, (1,1)}_{\wt\PP^1\times\Afuhs}.
\end{multline}
Comparing then to \eqref{eq:Fouriertau}, we get the assertion of Proposition 
\ref{prop:kBkB}.
\end{proof}

\subsection{Twistor conjugation and twistor sesquilinear pairings}\label{subsec:twconj}
We now work on the projective line $\PP^1$ equipped with two affine charts $\Omega_0$ and $\Omega_\infty$, and we denote by $\hb$ (\resp $\hb'$) the coordinate in the chart $\Omega_0$ (\resp $\Omega_\infty$) with $\hb'=1/\hb$.

\subsubsection*{Twistor conjugation}
In this setting, we denote by $c$ the conjugation functor considered in \S\ref{subsec:Fsesqui} and, if $\sigma:\PP^1\to c\PP^1$ or $c\PP^1\to\PP^1$ is the map\footnote{where we use the notation $c$ for the usual map sending a complex number to its conjugate, not to be confused with the conjugation functor or the conjugation morphism, also denoted by $c$.} $\hb\mto-1/c(\hb)$ or $c(\hb)\mto-1/\hb$, we denote by $\ov{\phantom{X}}$ the functor $\sigma^* c$. We call it the geometric or \emph{twistor conjugation} functor.

For instance, if $\cH$ is a $\cO_{\Omega_0}$-module, then $\ov\cH$ is a $\cO_{\Omega_\infty}$-module.

In the following, we denote by $\bS$ the circle $\module{\hb}=1$ and by $\cO_\bS$ the sheaf-theoretic restriction $\cO_{\Omega_0|\bS}$ (which can be identified with the sheaf of complex valued real-analytic functions on $\bS$).

\subsubsection*{Twistor sesquilinear pairing and objects of $\RTriples(\pt)$}
Let $\cH',\cH''$ be two $\cO_{\Omega_0}$-modules. A (twistor) sesquilinear pairing between these modules will be by definition a $\cO_\bS$-linear morphism
\[
C:\cH'_{|\bS}\ootimes_{\cO_\bS}\ov{\cH''_{|\bS}}\to\cO_\bS.
\]
In \cite[\S2.1.b]{Bibi01c}, we have denoted by $\RTriples(\pt)$ the category of such triples $(\cH',\cH'',C)$.

We say that such a triple is \emph{integrable} if $\cH',\cH''$ are equipped with a meromorphic connection having a pole of order at most $2$ at $\hb=0$ and no other pole (\ie if they are equipped with an action of $\hb^2\partial_\hb$) and if the sesquilinear pairing $C$ satisfies
\begin{equation}\label{eq:intC}
\hb\partial_\hb C(m',\ov{m''})= C(\hb\partial_\hb m',\ov{m''})-C(m',\ov{\hb\partial_\hb m''}),
\end{equation}
where the action of $\hb\partial_\hb$ on $\cH_{|\bS}$ is defined as that of $\hbm\cdot \hb^2\partial_\hb$, and the action of $\hb\partial_\hb$ on $\cO_\bS$ is the natural one. We denote by $\RdTriples(\pt)$ the category of integrable triples.

Denote by $\cL'\subset\cH'_{|\hb\neq0}$, $\cL''\subset\cH''_{|\hb\neq0}$ the local systems $\ker\hb^2\partial_\hb$. The local system attached to $\ov{\cH''}$ on $\hb\neq0$ is then $\sigma^{-1}c\cL''$. Let us notice that, when restricted to $\bS$, $\sigma$ is equal to $\iota$ (introduced in \S\ref{subsec:Fsesqui}) and that, when restricted to the local systems, the sesquilinear pairing $C$ takes values in the constant sheaf $\CC_\bS\subset\cO_\bS$ defined as $\ker\hb\partial_\hb$. Therefore, the sesquilinear pairing of an object of $\RdTriples(\pt)$ is determined by the $\CC$-linear morphism (its restriction to horizontal sections):
\begin{equation}\label{eq:twsesqui}
C:\cL'_{|\bS}\otimes_\CC\iota^{-1}c\cL''_{|\bS}\to\CC_\bS.
\end{equation}

\subsubsection*{Polarized twistor structures of weight~$0$}
Let $\cH',\cH''$ be two vector bundles (of the same rank) on $\Omega_0$. We say that the object $(\cH',\cH'',C)$ of $\RTriples(\pt)$ is a twistor structure of weight~$0$ if $C$ defines a gluing between $\cH^{\prime\vee}$ (dual bundle) and $\ov{\cH''}$ (in other words, $C$ is nondegenerate) and if the resulting vector bundle on $\PP^1$ is trivial (as we assume that the weight is $0$).

A \emph{polarization} is then an isomorphism (that we usually assume to be the identity) between $\cH'$ and $\cH''$ (\cf \cite[\S2.1.c]{Bibi01c}) such that, if we set $\cH=\cH'=\cH''$, the sesquilinear pairing $C:\cH_{|\bS}\otimes_{\cO_\bS}\ov{\cH_{|\bS}}\to\cO_\bS$ is \emph{Hermitian} (\ie $C(m,\ov\mu)=\ov{C(\mu,\ov m)}$ for local sections $m$ of $\cH_{|\bS}$ and $\mu$ of $\sigma^{-1}\cH_{|\bS}$) and \emph{positive definite}, \ie there exists a $\CC$-vector space $H\subset \Gamma(\Omega_0,\cH)$ such that
\begin{itemize}
\item
$\cH=\cO_{\Omega_0}\otimes_\CC H$,
\item
$C$ sends $H\otimes_\CC\ov H$ to $\CC\subset\Gamma(\bS,\cO_\bS)$ (hence induces a Hermitian form in the usual sense on $H$),
\item
and is positive definite as such.
\end{itemize}

\subsubsection*{Twistorization}

We will often use the following procedure, that we call \emph{twistorization}, which replaces a usual sesquilinear pairing as in \S\ref{subsec:Fsesqui} with a twistor sesquilinear pairing.

We consider the affine line with coordinate $\theta$. Let us assume that we are given a free $\CC[\theta]$-module $G_0$ of finite rank, equipped with a connection $\nabla$ having Poincar\'e rank one at the origin (\ie a pole of order two) and no other pole. We do not make any assumption on the behaviour of the connection at infinity, that will be lost anyway. Let us also assume that $G\defin\CC[\theta,\theta^{-1}]\otimes_{\CC[\theta]}G_0$ is equipped with a sesquilinear pairing $s:G^\an\otimes_{\cO_{\CC^*}}\ov{\iota^+G^\an}\to\cC^{\infty}_{\CC^*}$ as in \eqref{eq:Ftk}, compatible with the connection.

\begin{definition}[Twistorization]\label{def:twistorization}
The twistorization $(\cG,\cG,C)$ of the data $(G_0,\nabla,s)$ is the following object of $\RdTriples(\pt)$:
\begin{itemize}
\item
$\cG=G_0^\an$, equipped with $\nabla^\an$ (the analytization of $(G_0,\nabla)$), and the variable $\theta$ is renamed as $\hb$;
\item
the twistor sesquilinear pairing $C:\cG_{|\bS}\otimes_{\cO_\bS}\ov{\cG_{|\bS}}\to\cO_{\bS}$ (where $\ov{\cG}$ is the twistor conjugate of $\cG$) is obtained from \eqref{eq:twsesqui}, where we take for $\cL'=\cL''$ the local system $\ker\nabla^\an$ and for $C$ the restriction of $s$ to this local system, when restricted to \hbox{$\bS=\{|\hb|=1\}$}.
\end{itemize}
We equip this object with the sesquilinear duality $\cS=(\id_{\cH},\id_{\cH})$.
\end{definition}

\subsection{Fourier-Laplace transform of a filtered $\Clt$-module with a sesquilinear pairing}\label{subsec:Fourierfiltered}

In this paragraph we associate to any holonomic $\Clt$-module $M$ equipped with a good filtration $F_\bbullet M$ and a sesquilinear pairing $k:M\otimes_\CC\ov M\to\cS'(\Afu)$ an object $(G_0,\nabla,s)$ as in Definition~\ref{def:twistorization}, in order, through the twistorization, to get an object $(\cH,\cH,\wh C)$ of $\RdTriples(\pt)$. Let us note that $(G^\an,\nabla,s)$ has yet be obtained in \eqref{eq:Ftk}, and we are left to define $G_0$.

\subsubsection*{Saturation of lattice in a holonomic $\Clt$-module}
Let $M$ be a holonomic $\Clt$-module and let $L$ be a lattice of $M$, that is, $L\subset M$ is a $\CC[t]$-submodule of finite type and $M=\CC[\partial_t]\cdot L$ (notice that the generic rank of $L$ as a $\CC[t]$-module may be strictly smaller than that of $M$). We set $G\defin M[\ptm]=\CC[t]\langle\partial_t,\ptm\rangle\otimes_{\Clt}M$ (it is known that $G$ is also holonomic as a $\Clt$-module) and we denote by $\muloc:M\to G$ the natural morphism (the kernel and cokernel of which are isomorphic to powers of $\CC[t]$ with its natural structure of left $\Clt$-module). Let us set
\begin{equation}\label{eq:satdtm}
G_0^{(L)}=\sum_{j\geq0}\partial_t^{-j}\muloc(L).
\end{equation}
This is a $\CC[\ptm]$-submodule of $G$. Moreover, because of the relation $[t,\ptm]=(\ptm)^2$, it is naturally equipped with an action of $\CC[t]$. If $M$ has a regular singularity at infinity, then $G_0^{(L)}$ has finite type over $\CC[\ptm]$ (\cf \cite[Th\ptbl V.2.7]{Bibi00}).

For any $\ell\geq0$, we set $L_\ell=L+\cdots+\partial_t^\ell L$. We also have
\[
G_0^{(L)}=\sum_{j\geq0}\partial_t^{-j}\muloc(L_j)=\tbigcup_{j\geq0}\partial_t^{-j}\muloc(L_j).
\]
Moreover, for any $\ell\in\ZZ$, let us define $G_\ell^{(L)}$ as $\partial_t^\ell G_0^{(L)}\subset G$. Then we have, for any $\ell\geq0$,
\begin{equation}\label{eq:Gk}
G_\ell^{(L)}=\sum_{j\geq0}\partial_t^{-j}\muloc(L_\ell)=\sum_{j\geq0}\partial_t^{-j}\muloc(L_{\ell+j}),
\end{equation}
and therefore $G=\CC[\partial_t]\cdot G_0^{(L)}$ (and thus $G=\CC[\partial_t]\cdot G_\ell^{(L)}$ for any $\ell\in\ZZ$).

\subsubsection*{The case of a filtered $\Clt$-module}
Let us now assume that $(M,F_\bbullet M)$ is a $\Clt$-module equipped with a good filtration. Let $p_0\in\ZZ$. We say that $F_\bbullet M$ \emph{is generated by $F_{p_0}M$} if, for any $\ell\geq0$, we have $F_{p_0+\ell}M=F_{p_0}M+\cdots+\partial_t^\ell F_{p_0}M$. In other words, setting $L=F_{p_0}M$, $L$ is a lattice \emph{and} $F_{p_0+\ell}M=L_\ell$ for any $\ell\geq0$. We notice that the $\CC[\ptm]$-module $G_{-p_0}^{(F_{p_0})}=\partial_t^{-p_0}G_0^{(F_{p_0})}$ does not depend on the choice of the index $p_0$, provided that the generating assumption is satisfied. Indeed, setting $L'=F_{p_0+1}M=L_1$, we have $G_0^{(L')}=\partial_tG_0^{(L)}=G_1^{(L)}$ by \eqref{eq:Gk}, hence $G_{-(p_0+1)}^{(F_{p_0+1})}=G_{-p_0}^{(F_{p_0})}$. We thus set
\begin{equation}\label{eq:defG0F}
G_0^{(F_\bbullet)}=G_{-p_0}^{(F_{p_0})}\quad\text{for some (or any) index $p_0$ of generation.}
\end{equation}
If we also set $\theta=\partial_t^{-1}$, then $G_0^{(F_\bbullet)}$ is a free $\CC[\theta]$-module which satisfies $G=\CC[\theta,\theta^{-1}]\otimes_{\CC[\theta]}G_0^{(F_\bbullet)}$ and which is stable by the action of $\theta^2\partial_\theta=t$.

\begin{definition}[The basic correspondence]\label{def:basic}
If $(M,F_\bbullet M,k)$ consists of
\begin{itemize}
\item
a holonomic $\Clt$-module $M$ which is \emph{regular at $\infty$},
\item
a good filtration $F_\bbullet M$ of $M$,
\item
a sesquilinear pairing $k:M\otimes_\CC\ov M\to\cS'(\Afu)$,
\end{itemize}
we associate to such data an object $(\cH,\cH,\wh C)$ of $\RdTriples(\pt)$ by using the twistorization process of Definition \ref{def:twistorization} applied to the filtered Laplace transform of $(M,F_\bbullet M)$:
\begin{itemize}
\item
we set $\cH=G_0^{(F_\bbullet),\an}$ (the analytization of the object defined by \eqref{eq:defG0F}), by renaming $\hb$ the variable $\theta=\tau^{-1}=\ptm$, and we define the action of $\hb^2\partial_\hb$ as being that of~$t$,
\item
the sesquilinear pairing \eqref{eq:sesqui} induced by $F_tk$ is now regarded as a sesquilinear pairing \eqref{eq:twsesqui}, and therefore defines an integrable sesquilinear pairing $\wh C:\cH_{|\bS}\otimes_{\cO_\bS}\ov{\cH_{|\bS}}\to\cO_\bS$.
\end{itemize}
\end{definition}

\skippointrait
\begin{remarques}\label{rem:kBkB}
\begin{enumerate}
\item\label{rem:kBkB1}
In such a correspondence, if $k$ is Hermitian, then so is $\wh C$.
\item\label{rem:kBkB2}
If $M$ is assumed to have only regular singularities, then we may replace the datum of $k$ with that of the topological $k_B$. According to Proposition \ref{prop:kBkB} and Lemma \ref{lem:itwopi}, the sesquilinear pairing $\wh C$ is induced by $\itwopi\wh{k_B}$.
\end{enumerate}
\end{remarques}

\subsection{A criterion for the polarizability of $(\cH,\cH,\wh C)$}\label{subsec:crit}
In order to understand the basic construction above, it is useful to associate to $(M,F_\bbullet M,k)$ an object of $\RdTriples(\Afu)$ before taking Fourier-Laplace transform. We will first work algebraically in the coordinates $t$ and $\hb$.

\subsubsection*{The Rees module of a good filtration}
Denote by $R_FM$ the Rees module $\oplus_kF_kM\hb^k$, where $\hb$ is a new variable. We have $R_FM\subset \CC[\hb,\hbm]\otimes_\CC M$ and moreover $\CC[\hb,\hbm]\otimes_{\CC[\hb]}R_FM=\CC[\hb,\hbm]\otimes_\CC M$. This is a $\CC[t,\hb]\langle\partiall_t\rangle$-module\footnote{We denote by $\CC[t,\hb]\langle\partiall_t\rangle$ the Rees ring associated to the filtration of $\Clt$ by the order of operators, where we forget the grading; it is the free algebra generated by the polynomial algebras $\CC[t,\hb]$ and $\CC[\partiall_t]$ modulo the relation $[\partiall_t,t]=\hb$.} ($\partiall_t$ acts as $\hb\otimes\partial_t$). It is integrable, the $\hb^2\partial_\hb$-action being the natural one.

We define the conjugate object in a mixed sense: we use the standard conjugation with respect to the $t$-variable and the twistor one with respect to the $\hb$-variable. In particular, the conjugate $\ov{\CC[\hb,\hbm]\otimes_\CC M}$ is $\CC[\hb,\hbm]\otimes_\CC \ov M$ with the twistor-conjugate structure of $\CC[\hb,\hbm]$-module.

If $k$ is a sesquilinear pairing on $M$, then we extend it by $\CC[\hb,\hbm]$-linearity as
\begin{equation}\label{eq:Calg}
C:\big(\CC[\hb,\hbm]\otimes_\CC M\big)\otimes_{\CC[\hb,\hbm]}\big(\ov{\CC[\hb,\hbm]\otimes_\CC M}\big)\to \CC[\hb,\hbm]\otimes_\CC\cS'(\Afu).
\end{equation}
Clearly, $C$ satisfies the integrability condition like \eqref{eq:intC}.

\subsubsection*{Extension to $\PP^1$ and analytization}
Recall that we denote by $(*\infty)$ the effect of tensoring with $\cO_{\PP^1}(*\infty)$. In particular $\cD_{\PP^1}(*\infty)$ denotes the sheaf of analytic differential operators localized away from infinity. We will similarly consider the sheaves $R_F\cD_{\PP^1}(*\infty)$ (Rees sheaf of ring associated to the filtration of $\cD_{\PP^1}$ by the order of differential operators, localized away from infinity), its subsheaf $R_F\cO_{\PP^1}(*\infty)=\CC[\hb]\otimes_\CC\cO_{\PP^1}(*\infty)$ and $\cR_{\cP^1}(*\infty)$ (analytization of $R_F\cD_{\PP^1}(*\infty)$ with respect to $\hb$) as in~\cite{Bibi01c}, where we denote by $\cP^1$ the space $\PP^1\times\Omega_0$ with its analytic topology.

To $M$ one associates the sheaf $\ccM(*\infty)$ of $\cD_{\PP^1}(*\infty)$-modules, and to $R_FM$ one associates $R_F\ccM(*\infty)$. The $\cR_{\cP^1}(*\infty)$-module obtained from $R_F\ccM(*\infty)$ by analytization with respect to $\hb$ is now denoted by $\cM(*\infty)$: we have, by definition, $\cM(*\infty)=\cO_{\cP^1}(*\infty)\otimes_{R_F\cO_{\PP^1}(*\infty)}R_F\ccM(*\infty)$.

Let $\Db_{\PP^1}$ be the sheaf of distributions on $\PP^1$. Then $\cS'(\Afu)$ is nothing but the space of global sections of $\Db_{\PP^1}(*\infty)$. Formula \eqref{eq:Calg} can be sheafified to produce a sesquilinear pairing
\[
C:\cM(*\infty)_{|\bS}\otimes_{\cO_\bS}\ov{\cM(*\infty)_{|\bS}}\to \Dbh{\PP^1}(*\infty),
\]
where $\Dbh{\PP^1}$ denotes the sheaf of distributions on $\PP^1\times\bS$ which are continuous with respect to $\hb\in\bS$ (see \eg \cite[\S0.5]{Bibi01c}).

\begin{maintheorem}\label{th:main}
Let $M$ be a regular holonomic $\Clt$-module equipped with a good filtration $F_\bbullet M$ and a sesquilinear pairing $k:M\otimes_\CC\ov M\to \cS'(\Afu)$. Let us assume that the associated object $(\cM(*\infty),\cM(*\infty),C)$ of $\RdTriples(\PP^1)(*\infty)$ is the localization away from $\infty$ of a polarized regular twistor $\cD$-module of weight~$0$ on $\PP^1$ with polarization $(\id,\id)$ (\cf \cite[Def\ptbl 4.1.2 and 4.2.1]{Bibi01c}). Then the object $(\cH,\cH,\wh C)$ associated to $(M,F_\bbullet M,k)$ is an integrable polarized twistor structure of weight~$0$ with polarization $(\id,\id)$.
\end{maintheorem}

Let us express the result without using the language of twistors when
$p_0=0$. Denote by $\theta$ the variable $\tau^{-1}$. Then $G_0$
associated to $M_0=F_0M$ is a free $\CC[\theta]$-module of finite rank
$\mu$, equipped with an action of $\theta^2\partial_\theta$, induced
by the action of $t$. Denote by $\Delta$ the closed disc
$\{\module{\theta}\leq1\}$. Then, \emph{there exists a
$\cO(\Delta)$-basis $\omegag=(\omega_1,\dots,\omega_\mu)$ of
$\Gamma(\Delta,G_0^\an)$ such that, for any $i,j=1,\dots,\mu$, the
function $S^1\ni\tau\mto F_tk(\omega_i,\ov{\omega_j})(\tau)$ is
constant and equal to the Kronecker symbol $\delta_{ij}$} (where 
$\ov{\omega_j}$ is computed in the twistor sense).

\medskip
Let us notice that the basis $\omegag$ is in general not contained in $G_0$ (which is a natural $\CC[\theta]$-submodule of $\Gamma(\Delta,G_0^\an)$), \ie cannot be obtained by an algebraic base change from a basis of $G_0$. In other words, $G_0$ and $\CC[\theta]\cdot\omegag$ correspond to distinct algebraic extensions of $G_0^\an$. Indeed, on the one hand, the connection $\partial_\theta$ on $G_0$ has a regular singularity at $\theta=\infty$ (as $\wh M$ has a regular singularity at $\tau=0$). On the other hand, the connection $\partial_\theta$ on $\CC[\theta]\cdot\omegag$ has an irregular singularity at $\theta=\infty$ as soon as it has an irregular singularity at $\theta=0$.

\begin{exemples}\label{ex:basicsbis}
Let us explain the main theorem in the two elementary examples~\ref{ex:basics}.
\begin{enumerate}
\item\label{ex:basicsbis1}
In Example \ref{ex:basics}\eqref{ex:basics1}, we can take 
$F_0M=\CC\cdot m$, which generates a good filtration. We have $\wh 
M=\CC[\tau]\cE^{-c\tau}$, $G=\CC[\tau,\tau^{-1}]\cE^{-c\tau}$, and $G_0$ is 
the sub $\CC[\theta]$-module generated by $m$. We search for 
$\omega\in G_0^\an$ of the form $f(\theta)m$. We have $F_tk(m,\ov 
m)=e^{\ov{c\tau}-c\tau}=e^{\ov{c/\theta}-c/\theta}$. When restricted 
to $\module{\theta}=1$, this is written as $e^{\ov c\theta-c/\theta}$. 
Going to the twistor variable $\theta\mto\hb$, and using twistor conjugation, 
this is written as $f(\hb)\ov{f(\hb)}$ with $f(\hb)=e^{\ov c\hb}$. 
We can then choose $\omega=e^{-\ov c\theta}m$.
\item\label{ex:basicsbis2}
In Example \ref{ex:basics}\eqref{ex:basics2}, $m$ remains a generator 
of $G_0$ and, as $\module{\tau}=1$ on $S^1$, we can take 
$\omega=\sqrt{\sfrac{\Gamma(-\alpha)}{\Gamma(\alpha+1)}}\,m$.
\end{enumerate}
\end{exemples}

\section{Proof of the main theorem}

Let $\ccM$ be a regular holonomic $\cD_{\PP^1}$-module equipped with a good filtration $F_\bbullet\ccM$ and a sesquilinear pairing $k:\ccM\otimes_\CC\ov\ccM\to\Db_{\PP^1}$. All these data can be localized away from $\infty$ and, by taking global sections, we obtain $(M,F_\bbullet M,k)$ as in the main theorem.

\subsection{The Rees module and its Laplace transform}
We consider the Rees module $R_F\ccM$ associated to the filtration $F_\bbullet\ccM$, and its analytization (with respect to~$\hb$) that we denote by $\cM$. The conjugation is now taken in the usual sense with respect to the variable of $\PP^1$ and in the twistor sense with respect to the variable~$\hb$ as in \S\ref{subsec:crit} (\cf \cite[\S1.5.a]{Bibi01c}).

As in \S\ref{subsec:crit}, we construct a sesquilinear pairing $C:\cMS\otimes_{\cO_\bS}\ov{\cMS}\to\Dbh{\PP^1}$ from~$k$. Localizing away from $\infty$ gives the situation considered in \S\ref{subsec:crit}.

The assumption made in the main theorem is that there exists $(\ccM,F_\bbullet\ccM,k)$ such that the object $(\cM,\cM,C)$ (equipped with the isomorphism $\cS=(\id,\id)$) is a polarizable regular twistor $\cD$-module of weight~$0$.

Then, by Theorem~1 in \cite{Bibi04} (with its correction in the erratum to \cite{Bibi04}), the Fourier-Laplace transform $(\wh\cM,\wh\cM,\wh C)$ of $(\cM,\cM,C)$ with polarization $(\id,\id)$ is a polarizable regular twistor $\cD$-module of weight~$0$ on the Fourier plane with variable $\tau$, equipped with its analytic topology (that is, forgetting the behaviour at $\tau=\infty$)\footnote{More precisely, the proof in \cite{Bibi04} is given when the twistor object $(\cM,\cM,C)$ is simple and supported on $\PP^1$; the case when it is supported on a point is easy, as it reduces to Example \ref{ex:basicsbis}\eqref{ex:basicsbis1}.}. In particular, its fibre at $\tau=1$ is a polarizable twistor structure of weight~$0$.

The proof of the main theorem therefore reduces to the identification
of this fibre with the object constructed in
\S\ref{subsec:Fourierfiltered}.

\subsection{Laplace transform of the Rees module of a good filtration}
We keep notation of \oldS\S\ref{subsec:Fourierfiltered} and \ref{subsec:crit}. The Laplace transform $\wh{R_FM}$ of $R_FM$ is by definition the $\CC[\hb]$-module $R_FM$ equipped with the action of $\CC[\tau,\hb]\langle\partiall_\tau\rangle$ where $\tau$ acts as $\partiall_t$ and $\partiall_\tau$ as $-t$.

\begin{lemme}\label{lem:RFMhat}
The localized Laplace transform $\CC[\tau,\tau^{-1},\hb]\otimes_{\CC[\tau,\hb]}\wh{R_FM}$ with its natural $\CC[\tau,\tau^{-1},\hb]\langle\partiall_\tau\rangle$-structure is isomorphic to $\CC[\tau,\tau^{-1}]\otimes_\CC G_0^{(F_\bbullet)}$ equipped with the following structure:
\begin{itemize}
\item
the $\CC[\tau,\tau^{-1}]$-structure is the natural one,
\item
the multiplication by $\hb$ is given by $\hb\cdot (\tau^\ell\otimes g)=\tau^{\ell+1}\otimes(\ptm g)$, \ie $\hb\cdot{}=(\tau\otimes\ptm)\cdot{}$,
\item
the action of $\partiall_\tau$ is given by $\partiall_\tau(\tau^\ell\otimes g)=\tau^\ell\otimes[(\ell\ptm-t)g]$, \ie
\[
\partiall_\tau\cdot{}=\hb\cdot(\partial_\tau\otimes1)-1\otimes t=\tau\partial_\tau\otimes\ptm-1\otimes t.
\]
\end{itemize}
\end{lemme}

We see in particular that the fibre of $\wh{R_FM}$ at $\tau=1$ is nothing but $G_0^{(F_\bbullet)}$ with the $\CC[\hb]$ action defined by $\hb\cdot g=\ptm g$.

\begin{proof}
By definition, $\wh{R_FM}$ is included in $\CC[\hb,\hbm]\otimes_\CC M$, and the action of $\tau$ is induced by $\hb\otimes\partial_t$. By localization we thus have
\[
\CC[\tau,\tau^{-1}]\otimes_{\CC[\tau]}\big(\CC[\hb,\hbm]\otimes_\CC M\big)=\CC[\hb,\hbm]\otimes_\CC M[\ptm],
\]
where $\tau$ still acts as $\hb\otimes\partial_t$. The localized module $\CC[\tau,\tau^{-1},\hb]\otimes_{\CC[\tau,\hb]}\wh{R_FM}$ is therefore equal to the submodule of $\CC[\hb,\hbm]\otimes_\CC M[\ptm]$ generated by the $\tau^{-j}\muloc(\wh{R_FM})$ for $j\geq\nobreak0$. The coefficient of $\hb^\ell$ can be written as $\sum_{j\geq0}(\ptm)^j\muloc(F_{\ell+j})$. If $p_0$ is a generating index, this is nothing but $\partial_t^\ell G_{-p_0}^{F_{p_0}}$ for any $\ell\in\ZZ$ (for $\ell\geq p_0$ use \eqref{eq:Gk}, and for $\ell<p_0$, use also that, for any $k\in\ZZ$, we have $\muloc F_k\subset\ptm\muloc F_{k+1}$ because $\partial_tF_k\subset F_{k+1}$). The localized Rees module can now be written as $\oplus_{\ell\in\ZZ} (\partial_t^\ell G_0^{(F_\bbullet)}) \hb^\ell$, and the correspondence $\tau^\ell\otimes g\mto\partial_t^\ell g\hb^\ell$ induces an isomorphism $\CC[\tau,\tau^{-1}]\otimes_\CC G_0^{(F_\bbullet)}\to\oplus_{\ell\in\ZZ} (\partial_t^\ell G_0^{(F_\bbullet)})\hb^\ell$. One checks that the action of $\CC[\tau,\tau^{-1},\hb]\langle\partiall_\tau\rangle$ corresponds to that given in the lemma.
\end{proof}

\begin{remarque}[Integrability]\label{rem:IntRFM}
As $R_FM$ is naturally equipped with an action of $\hb^2\partial_\hb$ ($m_\ell\hb^\ell\mto \ell m_\ell\hb^{\ell+1}$), its Laplace transform $\wh{R_FM}$ is equipped with the twisted action
\[
m_\ell\hb^\ell\mto \ell m_\ell\hb^{\ell+1}+\partiall_ttm_\ell\hb^\ell=(\partial_tt+\ell)m_\ell\hb^{\ell+1};
\]
in order to understand this action, recall that $\wh{R_FM}$ can be regarded as the cokernel of the map
\[
e^{t\tau/\hb}\circ\partiall_t\circ e^{-t\tau/\hb}:\CC[\tau]\otimes_\CC R_FM\to\CC[\tau]\otimes_\CC R_FM;
\]
the action of $\partiall_\tau$, \resp $\hb^2\partial_\hb$, on this cokernel is induced by that of $e^{t\tau/\hb}\circ\partiall_\tau\circ e^{-t\tau/\hb}$, \resp $e^{t\tau/\hb}\circ\hb^2\partial_\hb\circ e^{-t\tau/\hb}$; the latter is written as $\hb^2\partial_\hb+t\tau$ and is translated by the formula above. The localized Laplace transform also has such an action. On the model $\CC[\tau,\tau^{-1}]\otimes_\CC G_0^{(F_\bbullet)}$, where the multiplication by $\hb$ is given by the action of $\tau\otimes\ptm$, the action of $\hb^2\partial_\hb$ is given by that of $\tau\otimes t$.

We thus see that the fibre of $\wh{R_FM}$ at $\tau=1$, as a $\CC[\hb]\langle\hb^2\partial_\hb\rangle$-module, is identified with $G_0^{(F_\bbullet)}$ with its natural structure of $\CC[\theta]\langle\theta^2\partial_\theta\rangle$-module (recall that we set $\theta=\ptm$, $\theta^2\partial_\theta=t$): $\hb$~acts as $\theta$ and $\hb^2\partial_\hb$ as $\theta^2\partial_\theta$.
\end{remarque}

At this stage, we have identified the fibre at $\tau=1$ of $\wh\cM$ (analytization of $\wh{R_FM}$) with its $\hb^2\partial_\hb$-action, with $G_0^{(F_\bbullet),\an}$ (analytization of $G_0^{(F_\bbullet)}$ with respect to $\theta$) with its $\theta^2\partial_\theta$-action.

\subsection{Identification of the sesquilinear pairings}

For any twistor $\cD$-module, the sesquilinear pairing is defined away from $\hb=0$. Let us begin thus by giving a more precise identification of $\wh{R_FM}$ on the domain $\hb\neq0$, $\tau\neq0$.

Let us localize with respect to $\hb$ the module considered in Lemma \ref{lem:RFMhat}. If we first localize $\wh{R_FM}$ with respect to $\hb$, we obtain the module $\CC[\hb,\hbm]\otimes_\CC M$. Then, localizing with respect to $\tau$ gives $\CC[\hb,\hbm]\otimes_\CC M[\ptm]$. We also have a description of this module as $\CC[\tau,\tau^{-1}]\otimes_\CC G$ if we set $G=M[\ptm]$. It comes equipped with an action $\CC[\tau,\tau^{-1},\hb,\hbm]\langle\partial_\tau,\partial_\hb\rangle$ by localization of that on $\wh{R_FM}$.

Let us denote $\theta=\tau^{-1}$ and set $\eta=\ptm$ acting on $G$, so that $t$ acts as $\eta^2\partial_\eta$. Let us give the explicit form of the action on $\CC[\tau,\tau^{-1}]\otimes_\CC G=\CC[\theta,\theta^{-1}]\otimes_\CC G$. The action of~$\hb$ is by $\theta^{-1}\otimes\eta$, that of $\partial_\hb$ is by $\theta\otimes\partial_\eta$. On the other hand, the action of $\theta$ is by $\theta\otimes1$ and that of $\partial_\theta$ is by $\partial_\theta\otimes1+\theta^{-1}\otimes\eta\partial_\eta= \partial_\theta\otimes1+\hb(1\otimes\partial_\eta)$ (using that $\partial_\theta=-\tau^2\partial_\tau$).

Consider the morphism $p:\CC^*\times\CC^*\to\CC^*$ defined by $(\theta,\hb)\mto\eta=\hb\theta$. Then the previous module $\CC[\theta,\theta^{-1}]\otimes_\CC G$ is nothing but the inverse image (in the sense of $\cD$-modules) $p^+G$, if $G$ is regarded as a $\CC[\eta,\eta^{-1}]\langle\partial_\eta\rangle$-module.

In particular, if $\wh\cV$ denotes the local system attached to $\wh V=G^\an$ on $\CC^*$, we see that the local system attached to $\wh{R_FM}$ on $\CC^*\times\CC^*$ is $p^{-1}\wh\cV$.

\medskip Let $u$ be a temperate distribution on $\Afu$ (coordinate
$t$). One considers its Fourier transform $\ccF_\hb u$ with parameter
$\hb\in\bS$ and kernel $\exp[-2i\im(t\tau/\hb)]\itwopi dt\wedge d\ov t$ 
(also written as $\exp[-2i\im(t\tau/\hb)]\frac1{2\pi i} \frac{dt}{\hb}\wedge 
\ov{\frac{dt}{\hb}}$). It belongs to the
space $\cS'(\Afuh\times\bS/\bS)$, \ie is a temperate distribution on
the product $\Afuh\times\bS$ which depends continuously on
$\hb\in\bS$.

Let $m,\mu\in M$ and $u=k(m,\ov\mu)\in\cS'(\Afu)$. Then, when
restricted to $\tau\neq0$, the distribution $\ccF_\hb u$ is $C^\infty$
with respect to $(\tau,\hb)$ (being locally part of a horizontal
section of an integrable connection). Working with the variable
$\theta=\tau^{-1}$, we see that $\ccF_\hb u$ is the inverse image
\emph{via} the map $p:(\theta,\hb)\mto\eta=\hb\theta$ of the usual
Fourier transform of the distribution $u$ restricted at
$\eta^{-1}\neq0$ (the kernel is $\exp -2i\im (t/\eta) \itwopi dt\wedge 
d\ov t$).

The sesquilinear pairing $k:M\otimes_\CC\ov M\to\cS'(\Afu)$ extends to a sesquilinear pairing (where conjugation is taken in the twistor sense with respect to the variable $\hb$)
\begin{equation}\label{eq:RFk}
\begin{aligned}
R_Fk:R_FM\otimes_{\CC[\hb,\hbm]}\ov{R_FM}&\to\CC[\hb,\hbm]\otimes_\CC\cS'(\Afu),\\
m_p\hb^p\otimes\ov{m_q\hb^q}&\mto\hb^p\ov\hb^qk(m_p,\ov{m_q})=(-1)^q\hb^{p-q}k(m_p,\ov{m_q}).
\end{aligned}
\end{equation}
By restricting \eqref{eq:RFk} to~$\bS$ we can define $\ccF_\hb R_Fk:\wh{R_FM}_{|\bS}\otimes_{\cO_\bS}\ov{\wh{R_FM}_{|\bS}}\to \cS'(\Afuh\times\bS/\bS)$ by composing $R_Fk$ with the Fourier transform $\ccF_\hb$. If one also restricts to $\tau\neq0$, one sees that such a sesquilinear pairing takes values in $C^\infty$ functions. Let us notice that $\ccF_\hb R_Fk$ was denoted $\wh{C}$ in \cite[\S8.2.c]{Bibi01c}, as the definition of direct image of a twistor $\cD$-module involves a factor $1/2\pi i$ in the pairing.

We wish to show that the restriction of $\ccF_\hb R_Fk$ to $\tau=1$ (that is, $\theta=1$) induces on $G_0^\an$ the pairing defined by the basic correspondence.

Recall that we set $S^1=\{\module{\tau}=1\}=\{\module{\theta}=1\}$ and $\bS=\{\module{\hb}=1\}$. The restriction to horizontal sections on $S^1\times\bS$ of the pairing $\ccF_\hb R_Fk$ is a sesquilinear pairing $(p^{-1}\wh\cV_{S^1\times\bS}\otimes_\CC\sigma^{-1}c(p^{-1}\wh\cV)_{S^1\times\bS}\to\CC$, as $p^{-1}\wh\cV$ is the sheaf of horizontal sections of $\wh{R_FM}$ on $\CC^*\times\CC^*$. Recall also that, on $\bS$, we have $\sigma(\hb)=-1/c(\hb)=-\hb=\iota(\hb)$.

If we use the involution $\iota:\eta\to-\eta$, we see that
$(\ccF_\hb R_Fk)_{|S^1\times\bS}$ is $p^{-1}$ of the sesquilinear pairing
$(F_tk)_B:\wh\cV_{|S^1}\otimes_\CC\iota^{-1}\ov{\wh\cV_{|S^1}}\to\CC_{S^1}$. 
The restriction to $\theta=1$ (that is, $\tau=1$) of $\ccF_\hb R_Fk$ thus
coincides with $(F_tk)_B$ at the level of horizontal sections, as was
to be proved.\qed

\subsection{Dilatation}\label{subsec:dilatation}
We will now identify the fibre at $\tau=\tau_o\neq0$ of the twistor object $(\wh{R_FM},\wh{R_FM},\ccF_\hb R_Fk)$ by a similar computation.

Let us fix $\lambda_o\in\CC^*$. We denote by $\mu_{\lambda_o}:\Afu\to\Afu$ the multiplication by $\lambda_o$. If $\varphi$ is a function on $\Afu$, we have $(\mu_{\lambda_o}^*\varphi)(t)=\varphi\circ\mu_{\lambda_o}(t)=\varphi(\lambda_ot)$.

If $F$ is a $\CC[t]$-module, we identify $\mu_{\lambda_o}^*F$ with the $\CC$-vector space $F$ where, for $m\in F$, the $\CC[t]$-action is given by $p(t)\cdot m=p(t/\lambda_o)m$. The fibre of $\mu_{\lambda_o}^*F$ at $t=1$ is identified with the fibre of $F$ at $t=\lambda_o$.

Similarly, if $M$ is a $\Clt$-module, the module $\mu_{\lambda_o}^+M$ is the $\CC$-vector space $M$ with the previous twisted action of $\CC[t]$ and $\partial_t\cdot m=\lambda_o\partial_tm$.

We then clearly have $\ov{\mu_{\lambda_o}^+M}=\mu_{\ov{\lambda_o}}^+\ov M$ and $\wh{\mu_{\lambda_o}^+M}=\mu_{1/\lambda_o}^+\wh M$.

If $u\in\cS'(\Afu)$, we define $\mu_{\lambda_o}^*u$ so that, for any $\varphi\in\cS(\Afu)$, one has $\langle \mu_{\lambda_o}^*u,\varphi\rangle=\module{\lambda_o}^{-2}\langle u,\mu_{1/\lambda_o}^*\varphi\rangle$. We have $\wh{\mu_{\lambda_o}^*u}=\module{\lambda_o}^{-2}\mu_{1/\lambda_o}^*\wh u$.

If $k:M'\otimes_\CC \ov{M''}\to\cS'(\Afu)$ is a sesquilinear pairing, we define $\mu_{\lambda_o}^*k$ so that, for any $m',m''\in M',M''$, we have $(\mu_{\lambda_o}^*k)(m',\ov m'')=\mu_{\lambda_o}^*\big(k(m',\ov m'')\big)$. Then $\mu_{\lambda_o}^*k$ is a sesquilinear pairing on $\mu_{\lambda_o}^+M'\otimes_\CC\ov{\mu_{\lambda_o}^+M''}$.

Given a filtered $\Clt$-module with sesquilinear pairing $(M,F_\bbullet M,k)$, we denote by $\mu_{\lambda_o}^+(M,F_\bbullet M,k)$ the object $(\mu_{\lambda_o}^+M,\mu_{\lambda_o}^*F_\bbullet M,\mu_{\lambda_o}^*k)$.

Similar formulas can be obtained for the dilatation of a Rees module and its Laplace transform. We then obtain:

\begin{lemme}
The fibre at $\tau=\tau_o\neq0$ of the twistor object $(\wh{R_FM},\wh{R_FM},\ccF_\hb R_Fk)$ is the object associated, through the basic correspondence of Definition \ref{def:basic}, to the triple $\mu_{1/\tau_o}^+(M,F_\bbullet M,k)$.\qed
\end{lemme}

\begin{remarque}\label{rem:dilat}
With this interpretation of $\mu_{1/\tau_o}^+(M,F_\bbullet M,k)$, Theorem~1 in \cite{Bibi04} can be restated by saying that the family $\mu_{1/\tau_o}^+(M,F_\bbullet M,k)$ ($\tau_o\in\CC$) corresponds, through the basic correspondence \ref{def:basic}, to a variation of polarized twistor structure of weight $0$ having a tame behaviour when $\tau_o\to0$.
\end{remarque}

\section{Variations of polarized~complex Hodge structure and polarized~twistor~$\cD$-modules}\label{sec:twVHSP}

We will associate to any variation of polarized complex Hodge structure\footnote{We restrict the study to weight $0$ for the sake of simplicity.} of weight~$0$ on $\PP^1\moins P$ (where $P$ is a finite set of points) a set of data $(M,F_\bbullet M,k)$ as in the basic correspondence of Def\ptbl \ref{def:basic}. In \S\ref{subsec:twVHSP} we will show that the assumptions of the Main Theorem are satisfied for these data. This will allow us to apply it to a variation of polarized complex Hodge structure and obtain Corollary \ref{cor:main}, which is the main result of this section.

The properties we want would basically follow from Schmid results \cite{Schmid73} (at least if we assume that the local monodromies of the variation are quasi-unipotent; for variations defined over $\RR$, see \cite[\S11]{Zucker79}). Nevertheless, we will directly construct the twistor $\cD$-module, using the general results of \cite{Simpson90}, as translated in terms of twistor $\cD$-modules in \cite[Chap\ptbl5]{Bibi01c}. Our objective is to make clear the characterization of those polarized twistor $\cD$-modules which come from a variation of Hodge structure. The study of smooth objects is made in \S\ref{subsec:smoothtwVHSP} and their extension to $\PP^1$ in \S\ref{subsec:twVHSP}.

The variation of polarized complex Hodge structure we start with is a set of data defined on $\PP^{1\an}\moins P$ (\cf \S\ref{subsec:varHodge}). We want to extend these data to $\PP^1$. The problem is local near each puncture in $P$, so we work locally analytically near each puncture. We denote by $X$ the disc of radius $1$ centered at the origin in the complex plane with coordinate $x$ and we denote by $X^*$ the punctured disc $X\moins\{0\}$.

\subsection{Variation of polarized complex Hodge structure on $X^*$}\label{subsec:varHodge}
We consider on $X^*$ a variation of complex Hodge structure of weight~$0$, which is polarized. It consists in giving a $C^\infty$ vector bundle $H$ on $X^*$ equipped with a flat connection $D$, a decomposition $H=\oplus_{p\in\ZZ}H^p$ ($H^p$ is usually written as $H^{p,-p}$ as the weight is $0$) and a Hermitian metric $h$ on $H$, satisfying the following properties:
\begin{itemize}
\item
the decomposition is orthogonal with respect to $h$ and the nondegenerate Hermitian form $k=\oplus_p(-1)^ph_{|H^p}$ is $D$-flat,
\item
(Griffiths' transversality)
\begin{equation}\label{eq:griffiths}
\begin{split}
D'(H^p)&\subset(H^p\oplus H^{p-1})\otimes_{\cO_{X^*}}\Omega^1_{X^*}\\
D''(H^p)&\subset(H^p\oplus H^{p+1})\otimes_{\cO_{\ov X^*}}\Omega^1_{\ov X^*}.
\end{split}
\end{equation}
\end{itemize}

We define the (increasing) Hodge filtration $F_\bbullet H$ as
\[
F_pH=\oplus_{q\geq -p}H^q,
\]
so that $D'F_pH\subset F_{p+1}H\otimes_{\cO_{X^*}}\Omega^1_{X^*}$.

We denote by $(V,\nabla)$ the holomorphic bundle with connection $(\ker D'',D')$ and we set $F_pV=F_pH\cap V$. We have $\nabla F_pV\subset F_{p+1}V\otimes_{\cO_{X^*}}\Omega^1_{X^*}$.

The triple $(H,D,h)$ is \emph{harmonic} (\cf \cite{Simpson90}), the metric connection $D_h$ (\resp the Higgs field~$\theta$) is obtained by composing $D$ with the projection on the first (\resp second) factor in \eqref{eq:griffiths}. In particular, $D_h$ respects the decomposition.

\begin{remarque}[Twist and shift]\label{rem:twsh}
Let $w$ be some integer. A variation of polarized complex Hodge structure 
of weight $w$ consists of the same data $H=\oplus_p H^p$ and the 
Hermitian metric $h$, where one now defines $k=i^w(-1)^ph$ on $H^p$ (which is written as 
$H^{p,w-p}$), satisfying the same properties. Going from weight $0$ 
to weight $w$ consists only in multiplying $k$ by $i^w$.

If one shifts the filtration $F_\bbullet$ or the decomposition 
$\oplus_pH^p$ by setting $F[-\ell]_p=F_{p+\ell}$ or 
$H[-\ell]=\oplus_pH^{p-\ell}$ ($\ell\in\ZZ$), the pairing $k$ has to be replaced 
with $k'=(-1)^\ell k$.
\end{remarque}

\subsection{The $\cD$-module associated to a variation of polarized complex Hodge structure}\label{subsec:Dmodvar}

We first extend the holomorphic vector bundle with connection $(V,\nabla)$ as a meromorphic vector bundle with regular connection on $X$. We denote by \hbox{$j:X^*\hto X$} the open inclusion.

\pagebreak[3]
\skippointrait
\begin{theoreme}[\cf \cite{Schmid73}]\label{th:schmid}
\begin{enumerate}
\item\label{th:schmid1}
The $\cO_X[1/x]$-module $(j_*V)^{\rmod}\subset j_*V$ of local sections of $j_*V$ on $X$, the $h$-norm of which has moderate growth near the origin, is free of rank $\rk V$. We denote it by $\ccM[1/x]$.
\item\label{th:schmid2}
The connection $\nabla$ is meromorphic on $\ccM[1/x]$ and has a regular singularity at the origin.
\item\label{th:schmid3}
Each term $(j_*V)^b$ ($b\in\RR$) of the parabolic filtration of $h$ at the puncture---defined as the $\cO_X$-module of local sections $v$ of $j_*V$ such that, for any $\epsilon>0$, $\module{x}^{-b+\epsilon}\norme{v}_h$ is bounded near the origin--- is a locally free $\cO_X$-module of rank $\rk V$.
\item\label{th:schmid4}
The connection $\nabla$ is logarithmic on each $(j_*V)^b$. On $(j_*V)^b/(j_*V)^{>b}$, its residue has $b$ as unique eigenvalue. In particular, the filtration $(j_*V)^\bbullet$ of $(j_*V)^{\rmod}$ coincides with the (decreasing) Malgrange-Kashiwara filtration $V^b(\ccM[1/x])$\footnote{\label{footnoteMK}Recall that, as $x$ is invertible on $\ccM[1/x]$, the Malgrange-Kashiwara filtration of $\ccM[1/x]$ is the unique filtration by locally free $\cO_X$-modules of finite rank such that the connection $\nabla$ on $V^b$ has a logarithmic singularity at $0$ with residue having $b$ as unique eigenvalue.}.
\end{enumerate}
\end{theoreme}

The parabolic filtration $(j_*V)^\bbullet$ is decreasing and we will usually consider the associated increasing filtration $(j_*V)_a=(j_*V)^{-a-1}$ for $a\in\RR$.

\begin{proof}[Sketch of proof]
One reduces to variations defined over $\RR$ in a standard way. Then the result is essentially proved in \cite{Schmid73} (\cf also \cite[\S11]{Zucker79}).
\end{proof}

\begin{remarque}
The basic result of Schmid that the $h$-norm of $D$-horizontal sections has moderate growth near the origin is equivalent to saying that the harmonic bundle $(H,D,h)$ is \emph{tame} in the sense of \cite{Simpson90}. Theorem \ref{th:schmid} is proved in this more general setting of tame harmonic bundles in \loccit We will use the more general version stated at the level of twistor $\cD$-modules in \S\ref{subsec:twVHSP}.
\end{remarque}

The $\cO_X[1/x]$-module with connection $(\ccM[1/x],\nabla)$ is thus a $\cD_X$-module with regular singularity at~$0$. Denote by $\ccM$ its minimal extension at the origin: by definition, this is the $\cD_X$-submodule generated by $(j_*V)_{<0}$ in $\ccM[1/x]$. For $a<0$, the filtration $(j_*V)_a$ of $\ccM$ is nothing but the $V$-filtration of Malgrange-Kashiwara at the origin, that we denote\footnote{The letter $V$ here should not be confused with the previous notation, corresponding to the holomorphic bundle; it is the standard notation for the Malgrange-Kashiwara filtration.} by $V_a\ccM$.

\subsection{Extending the sesquilinear pairings}\label{subsec:extsesqui}
By definition (and by Cauchy-Schwarz), the metric $h$ extends as a $\cO_X\otimes_\CC\cO_{\ov X}$-linear pairing $V_{<0}\ccM\otimes_\CC\ov{V_{<0}\ccM}\to L^1_{\loc,X}(\vol)$, where $\vol$ is the Euclidean volume form on $X$. Unless $h$ is flat, this sesquilinear form, regarded as taking values in the sheaf $\Db_X$ of distributions, is not $\cD_X\otimes_\CC\cD_{\ov X}$-linear.

On the other hand, the sesquilinear form $k$ also extends as a sesquilinear pairing $V_{<0}\ccM\otimes_\CC\ov{V_{<0}\ccM}\to L^1_{\loc,X}(\vol)$: this is seen, using Schmid's results, by considering the matrix of the base change between a horizontal basis of $V$, where $k$ is constant by flatness, and a basis of $V_{<0}\ccM$; this will be also recovered in \S\ref{subsec:twVHSP} where we will also obtain:

\begin{proposition}
The sesquilinear form $k$ is $V_0\cD_X\otimes_\CC V_0\cD_{\ov X}$-linear and extends in a unique way as a $\cD_X\otimes_\CC\cD_{\ov X}$-linear form
\[
k:\ccM\otimes_\CC\ov{\ccM}\to\Db_X.
\]
\end{proposition}

\begin{remarque}
To be more precise, the pairing $k$ takes values in the subsheaf of \emph{regular holonomic distributions} (\cf \cite{Kashiwara86,Bjork93} for such a notion, which will not be used here).
\end{remarque}

\subsection{Extending the Hodge filtration}\label{subsec:extendingHodge}
We wish to define a good filtration $F_\bbullet \ccM$ on~$\ccM$, starting from the Hodge filtration $F_\bbullet V$. We follow \cite[\S3.2]{MSaito86}. We first set, for any $p\in\ZZ$,
\begin{equation}\label{eq:FVneg}
F_pV_{<0}\ccM=j_*F_pV\cap V_{<0}\ccM.
\end{equation}
This is clearly a $\cO_X$-module (recall that $V_{<0}\ccM=(j_*V)_{<0}$). In other words, a local section $v$ of $j_*V$ on $X$ is in $F_pV_{<0}\ccM$ if and only if
\begin{itemize}
\item
$v$ is a local section of $j_*F_pV$,
\item
$\lim_{x\to0}\module{x}\cdot\norme{v}_h=0$.
\end{itemize}
Let us notice that $F_pV_{<0}\ccM=V_{<0}\ccM$ for $p\gg0$ and that $F_pV_{<0}\ccM=0$ for $p\ll0$. We next define, for any $p\in\ZZ$,
\begin{equation}\label{eq:FccM}
F_p\ccM=\sum_{j\geq0}\partial_x^j\cdot F_{p-j}V_{<0}\ccM\subset \ccM.
\end{equation}
This is clearly a $\cO_X$-module, we have $\partial_xF_p\ccM\subset F_{p+1}\ccM$ for any $p\in\ZZ$, and $F_p\ccM=0$ for $p\ll0$. Let us also notice that $F_p\ccM\cap V_{<0}\ccM=F_pV_{<0}\ccM$, as we have $F_p\ccM\subset j_*F_pV$, and that $j^{-1}F_p\ccM=F_pV$.

For any $a<0$ we have
\begin{equation}\label{eq:FpVa}
F_pV_{<0}\ccM\cap V_a\ccM=F_p\ccM\cap V_a\ccM=j_*F_pV\cap V_a\ccM.
\end{equation}
Indeed, the inclusions $\subset$ are clear. On the other hand, we have $j_*F_pV\cap V_a\ccM=(j_*F_pV\cap V_a\ccM)\cap V_{<0}\ccM=F_pV_{<0}\ccM\cap V_a\ccM$. For such an $a$, there is no ambiguity to denote by $F_pV_a\ccM$ any of the expression in \eqref{eq:FpVa}.

\skippointrait
\begin{proposition}\label{prop:Fcoh}
\begin{enumerate}
\item\label{prop:Fcoh1}
For any $p\in\ZZ$, $F_p\ccM$ is $\cO_X$-coherent.
\item\label{prop:Fcoh2}
The filtered $\cD_X$-module $(\ccM,F_\bbullet\ccM)$ is strictly specializable at $x=0$, that is,
\begin{enumerate}
\item\label{prop:Fcoh2a}
for any $a<0$ and any $p\in\ZZ$, we have $x\cdot F_pV_a\ccM=F_pV_{a-1}\ccM$,
\item\label{prop:Fcoh2b}
for any $a\geq-1$ and any $p$, we have $\partial_x\cdot F_p\gr_a^V\ccM=F_{p+1}\gr_{a+1}^V\ccM$, where we set, as usual, $F_p\gr_a^V\ccM=(F_p\ccM\cap V_a\ccM)/(F_p\ccM\cap V_{<a}\ccM)$.
\end{enumerate}
\end{enumerate}
\end{proposition}

It follows from \eqref{prop:Fcoh1} that $F_p\ccM$ is (locally) free as a $\cO_X$-module, as it has no $\cO_X$-torsion (being contained in $\ccM[1/x]$).

\begin{proof}[Sketch of proof]
The $\cO_X$-coherence of $F_pV_{<0}\ccM$ is the main point. It can be obtained from Schmid's Nilpotent Orbit Theorem \cite{Schmid73}, but we will recover it in \S\ref{subsec:twVHSP}. The coherence of $F_p\ccM$ follows, hence \eqref{prop:Fcoh1}. For \eqref{prop:Fcoh2}, argue as in \cite[Prop\ptbl 3.2.2]{MSaito86}.
\end{proof}

\subsection{The smooth polarized twistor structure associated to a variation of polarized~complex Hodge structure}\label{subsec:smoothtwVHSP}

Let $(H,D,k)$ be a variation of polarized complex Hodge structure of weight~$0$ on a Riemann surface $Y$ (we will take $Y=\PP^1\moins P$ or $Y=X^*$, the punctured disc). Let $(V,\nabla)$ be the corresponding holomorphic bundle with holomorphic connection and $F_\bbullet V$ its increasing Hodge filtration, as in \S\ref{subsec:varHodge}.

We associate with $(H,D,k)$ as above the triple $\cT=(R_FV,R_FV,R_Fk)$, where $R_FV$ is the Rees module $\oplus_{k\in\ZZ}\hb^kF_kV$ and $R_Fk$ is the sesquilinear pairing obtained from $k$ by twistor sesquilinearity over $\CC[\hb,\hbm]$ as in \eqref{eq:RFk}:
\[
\textstyle R_Fk\big(\sum_p\hb^pv_p,\ov{\sum_q\hb^qw_q}\big)\defin\dpl \sum_{p,q}\hb^p\ov\hb^qk(v_p,\ov{w_q}).
\]

\begin{lemme}
The triple $\cT=(R_FV,R_FV,R_Fk)$ equipped with the polarization $\cS=(\id,\id)$ is a smooth polarized twistor structure of weight~$0$.
\end{lemme}

\begin{proof}
It is enough to show that each fibre of $\cT$ over $Y$ is a polarized twistor structure of weight~$0$; we can therefore assume that $Y$ is a point, so that there is no difference between $V$ and~$H$. We are reduced to finding a $\CC[\hb]$-basis $\epsilong$ of $R_FV$ which is orthonormal for $R_Fk$ (\cf \cite[Remark~2.2.3]{Bibi01c}).

For any $p$, let $\epsilong_p$ be a $h$-orthonormal basis of $H_p$. Then $(\hb^p\epsilong_p)_{p\in\ZZ}$ is the desired basis of $R_FV$: that it is a basis is clear, and
\[
R_Fk(\hb^p\varepsilon_{p,i},\ov{\hb^q\varepsilon_{q,j}})=\hb^p\ov\hb^q k(\varepsilon_{p,i},\ov{\varepsilon_{q,j}});
\]
this expression vanishes unless $p=q$, hence is equal to $(-1)^pk(\varepsilon_{p,i},\ov{\varepsilon_{p,j}})=h(\varepsilon_{p,i},\ov{\varepsilon_{p,j}})=\delta_{i,j}$, if $\delta$ is the Kronecker symbol.
\end{proof}

\begin{remarque}[Integrability]
It is easily seen that this object of $\RTriples(Y)$ is integrable: indeed, $R_FV$ is integrable (\cf Remark \ref{rem:IntRFM}) and $R_Fk$ satisfies the integrability condition
\[
\hb\partial_\hb R_Fk(u,\ov v)=R_Fk(\hb\partial_\hb u,\ov v)- R_Fk(u,\ov{\hb\partial_\hb v}).
\]
\end{remarque}

\subsection{Characterization of polarized smooth twistor structures coming from variations of polarized complex Hodge structure}\label{subsec:charact}

Consider the analytization of the object $(R_FV,R_FV,R_Fk)$ constructed in the previous section, \ie tensor it with $\cO_\cY$ over $\cO_Y[\hb]$. It takes the form $(\cH',\cH',C)$, where $\cH'$ is a locally free $\cO_\cY$-module of finite rank. Is it possible to recover $(R_FV,R_FV,R_Fk)$ from its analytization, and how to do so?

\begin{proposition}\label{prop:smoothgraded}
Let $\cT=(\cH',\cH',C)$ be a smooth object of $\RTriples(Y)$. It is the analytization of a triple $(R_FV,R_FV,R_Fk)$ if and only if it satisfies the following supplementary properties:
\begin{enumerate}
\item\label{prop:smoothgraded1}
$\cT$ is integrable,
\item\label{prop:smoothgraded2}
$\cH'$ is relatively logarithmic, \ie stable under $\hb\partial_\hb$ (and not only under $\hb^2\partial_\hb$),
\item\label{prop:smoothgraded3}
the monodromy of the flat connection on $\cH'_{|\hb\neq0}$ around $\hb=0$ is the identity.
\end{enumerate}
\end{proposition}

\begin{proof}
The conditions are necessary: this clearly follows from the definition of the $\hb^2\partial_\hb$-action on $R_FV$ for \ref{prop:smoothgraded}\eqref{prop:smoothgraded1} and \eqref{prop:smoothgraded2}; localizing along $\hb=0$ (\ie tensoring with $\cO_Y[\hb,\hbm]$) changes $R_FV$ to $\CC[\hb,\hbm]\otimes_\CC V$ and the $\partial_\hb$-action is trivial on $1\otimes V$. This gives \eqref{prop:smoothgraded3}.

Let us now consider $\cT$ satisfying Properties \ref{prop:smoothgraded}\eqref{prop:smoothgraded1}--\eqref{prop:smoothgraded3}. We argue in four steps:
\begin{enumeratea}
\item\label{enum:a}
We show that, locally on $Y$, there exists a basis of $\cH'_{|\hb=0}$ as a $\cO_{\cY|\hb=0}$-module (where $_{|\hb=0}$ means the sheaf-theoretical restriction) such that the matrix of the connection in this basis is relatively logarithmic and takes the form $Dd\hb/\hb+B(y,\hb)dy/\hb$, where $D$ is a diagonal matrix with integral entries and $B$ is holomorphic.
\item\label{enum:b}
If $G_k$ denotes the free $\cO_Y$-module generated by the part of the previous basis corresponding with the eigenvalue $k$ of $D$, the $\cO_Y$-module $\cF_k\defin\ker[(\hb\partial_\hb-k):\cH'_{|\hb=0}\to\cH'_{|\hb=0}]$ is equal to $\bigoplus_{\ell\in\NN}G_{k-\ell}\hb^\ell$, hence is locally free. Moreover, the integrability of the connection implies that $\hb\partial_y\cF_k\subset \cF_{k+1}$. Setting $V=\ker[\hb\partial_\hb:\cH'[\hbm]_{|\hb=0}\to\cH'[\hbm]_{|\hb=0}$, we see that $V=\bigoplus_{\ell\in\ZZ}G_{\ell}\hb^{-\ell}$ is a locally free $\cO_Y$-module filtered by the $F_kV\defin \hb^{-k}\cF_k=\bigoplus_{\ell\leq k}G_\ell\hb^{-\ell}$, and that this filtration satisfies Griffiths transversality $\partial_yF_k\subset F_{k+1}$. Let us also notice that the filtration $F_\bbullet V$ is the filtration by the order of the pole when $V$ is considered as a subsheaf of $\cH'[\hbm]_{|\hb=0}$, that is, $F_kV=V\cap\hb^{-k}\cH'_{|\hb=0}$.
\item\label{enum:c}
The morphism $\cO_{\cY|\hb=0}\otimes_{\cO_Y}R_FV\to\cH'_{|\hb=0}$ sending $F_k\hb^k$ to $\cF_k$ is an isomorphism of germs of bundles with meromorphic connection. It extends in a unique way by horizontality as a morphism of bundles with connections $\cO_{\cY}\otimes_{\cO_Y}R_FV\to\cH'$. This morphism identifies $(V,\nabla)$, defined as $\ker\hb\partial_\hb$ as above, with $\cH'/(\hb-1)\cH'$ equipped with its natural connection.
\item\label{enum:d}
The integrability property for $C$ shows that, when restricted to $\cF_{k|\bS}\otimes_{\cO_\bS}\ov{\cF_{\ell|\bS}}$, $C$ is homogeneous of degree $(k-\ell)$ with respect to $\hb\partial_\hb$. Therefore, $C$ takes the form $R_Fk$ for some sesquilinear pairing $k$ on $V$.
\end{enumeratea}

Let us indicate the proof of Step \eqref{enum:a}. This is a particular case of the Levelt normal form with parameter. Take a local coordinate $y$ on $Y$ and choose a local basis of $\cH'_{|\hb=0}$. The matrix of the connection in this basis can be written as
\[
\wt A(y,\hb)\,\frac{d\hb}{\hb}+\wt B(y,\hb)\,\frac{dy}{\hb},
\]
where $\wt A$ and $\wt B$ are holomorphic. As the monodromy relative to $\hb$ is unipotent (being the identity), the characteristic polynomial of $\wt A(0,y)$ is constant and its roots are integers. Therefore, one can assume that $\wt A(0,y)=D+N(0,y)$, where $D$ is diagonal with integral eigenvalues and $N(0,y)$ is strictly upper triangular and commutes with~$D$.

Arguing as for the construction of the Levelt normal form (see \eg \cite[Exer\ptbl II.2.20]{Bibi00}), it is then possible to find a finite number of (possibly nonzero) holomorphic matrices $A_j(y)$ such that $[D,A_j]=-jA_j$ ($j\in\NN^*$) and a formal series
\[
\wh P(y,\hb)=\id+\hb P_1(y)+\cdots
\]
where $P_k(y)$ are holomorphic on a fixed neighbourhood of $y=0$ such that
\[
\hb\frac{\partial\wh P}{\partial\hb}(y,\hb)=\wh P(y,\hb)\cdot\big(D+N(0,y)+\textstyle\sum_j\hb^jA_j(y)\big)-\wt A(y,\hb)\wh P(y,\hb).
\]
After the formal meromorphic base change with matrix $\wh P(y,\hb)\hb^{-D}$, the matrix of $\hb\partial_\hb$ is $N(0,y)+\sum_jA_j(y)$ (see \eg \loccit) and therefore, restricting to curves $y=y_o$ and applying the classical theory of differential equations with regular singularities in dimension one, the monodromy around $\hb=0$ is $\exp-2\pi i(N(0,y)+\sum_jA_j(y))$. The assumption that the monodromy around $\hb=0$ is the identity is then equivalent to the vanishing of each term in the sum.

Applying now the formal base change with matrix $\wh P(y,\hb)$ instead of $\wh P(y,\hb)\hb^{-D}$, we find that in the new (formal) basis, the matrix of the connection is written as
\[
D\,\frac{d\hb}{\hb}+\wh B(y,\hb)\,\frac{dy}{\hb},
\]
where $\wh B= \wh P^{-1}\wt B\wh P+\hb\wh P^{-1}\partial\wh P/\partial y$ is a formal series in $\hb$ with holomorphic coefficients in $y$. In particular, we have $\wh B(0,y)=\wt B(0,y)$ and the integrability condition implies $\hb\partial\wh B(y,\hb)/\partial\hb=[\wh B,D]+\wh B$. Expanding this equality with respect to powers of $\hb$ shows that $\wh B$ is a polynomial in $\hb$ with holomorphic coefficients in $y$. The new matrix has the desired form.

Last, the formal matrix $\wh P$ is a horizontal section of a holomorphic bundle with meromorphic connection having regular singularities along $\hb=0$ (in order to justify this statement, let us remark that is so for its restriction to curves $y=y_o$; we then apply the regularity criterion in \cite{Deligne70}). It is therefore convergent.
\end{proof}

\subsection{The polarized twistor $\cD$-module associated to a variation of polarized~complex Hodge structure}\label{subsec:twVHSP}

Let $(H,D,k)$ be a variation of polarized complex Hodge structure of weight~$0$ on $\PP^1\moins P$. By the results of \oldS\S\ref{subsec:Dmodvar}, \ref{subsec:extsesqui} and \ref{subsec:extendingHodge} applied to the neighbourhood of each point of $P$, we associate to it a filtered $\cD$-module $(\ccM,F_\bbullet\ccM)$ on~$\PP^1$ equipped with a Hermitian pairing $k:\ccM\otimes_\CC\ov{\ccM}\to\Db_{\PP^1}$.

The main result of this section is:
\begin{proposition}\label{prop:twHodge}
The object of $\RTriples(\PP^1)$ associated to $(\ccM,F_\bbullet\ccM,k)$ is a polarized regular twistor $\cD$-module of weight~$0$ on $\PP^1$.
\end{proposition}

If $(M,F_\bbullet M)$ denotes the (global sections of) the localization away from $\infty$ of $(\ccM,F_\bbullet\ccM)$ and if $k$ denotes the associated Hermitian pairing with values in $\cS'(\Afu)$, it follows that the Main Theorem of \S\ref{subsec:crit} applies to $(M,F_\bbullet M,k)$:

\begin{corollaire}\label{cor:main}
Let $(H,D,k)$ be a variation of polarized complex Hodge structure of weight~$0$ on $\PP^1\moins P$, let $(M,F_\bbullet M)$ be the corresponding filtered $\Clt$-module and still denote by $k$ the extension of $k$ as a sesquilinear pairing $M\otimes_\CC\ov M\to\cS'(\Afu)$. Then the associated object $(\cH,\cH,\wh C)$ through the correspondence \ref{def:basic} is an integrable polarized twistor structure of weight~$0$ with polarization $(\id,\id)$.\qed
\end{corollaire}

In order to prove Proposition \ref{prop:twHodge}, we will directly construct a twistor $\cD$-module extending the one attached to the variation on $\PP^1\moins P$, and show that this object takes the form $(R_F\ccM,R_F\ccM,R_Fk)$ for some filtration $F$ on $\ccM$ and sesquilinear pairing $k$. We will then show that the filtration $F$ and the sesquilinear pairing $k$ coincide with those defined in \oldS\S\ref{subsec:Dmodvar}, \ref{subsec:extsesqui} and \ref{subsec:extendingHodge}. This will give in particular the finiteness results obtained there using Schmid's results.

\subsubsection*{Extending $(\cH',\cH',C)$}
Denote by $(\cH',\cH',C)$ the analytic triple attached to $(H,D,h)$ as in \S\ref{subsec:charact}, with polarization $(\id,\id)$. We have yet seen that the results of Schmid imply that the harmonic metric $h$ is tame near the punctures. In \cite[Cor\ptbl5.3.1]{Bibi01c} we have constructed, as a consequence of the results of \cite{Simpson90} and \cite{Biquard97}, an object $\cT=(\cM,\cM,C)$ which is a polarized twistor $\cD$-module if we take $(\id,\id)$ as the polarization, such that it restricts to $(\cH',\cH',C)$ on $\PP^1\moins P$. More precisely, we have also defined the extension $V_{<0}\wt\cM\subset j_*\cH'$, which is a locally free $\cO_{\cP^1}$-module (where $\cP^1=\PP^1\times\Omega_0$ and $\Omega_0$ is defined in \S\ref{subsec:twconj}), where $j:\PP^1\moins P\hto\PP^1$ denotes the open inclusion (and also the same inclusion after the product with $\Omega_0$). Let us notice that, as the eigenvalues of the local monodromies of $(H,D)$ near the punctures have modulus equal to one (\cf \cite{Schmid73}), the $V$-filtration constructed in \cite{Bibi01c} is defined globally with respect to $\hb$ and not only locally near each $\hb_o$. We will show that $(\cM,\cM,C)$ is the analytization of some $(R_F\ccM,R_F\ccM,R_Fk)$.

\subsubsection*{Extending $R_FV$}
Let us first show the existence of $R_FV_{<0}\ccM$.

\begin{lemme}
The $\hb\partial_\hb$ action on $\cH'$ extends to $V_{<0}\wt\cM$.
\end{lemme}

\begin{proof}
Let us work near a puncture, with local coordinate $x$. Recall (\cf \cite[Cor\ptbl5.3.1]{Bibi01c}) that $V_{<0}\wt\cM$ is defined as the subsheaf of $j_*\cH'$, the germ at $(0,\hb_o)$ of which consists of local sections $m$ such that $\lim_{x\to0}\module{x}\cdot\norme{m}_{\pi^*h}=0$, uniformly for $\hb$ near $\hb_o$.

Let $X$ denote a small disc centered at $x=0$, $X^*$ the punctured disc and $\ov X$ the closure of $X$. Choose an orthonormal basis $(\epsilon_{p,j})$ of $H$ on $\ov X^*$, which is adapted to the decomposition $H=\oplus_pH_p$, and write $m=\sum_{p,j}f_{p,j}(x,\hb)\hb^p\epsilon_{p,j}$, where each $f_{p,j}$ is $C^\infty$ with respect to $x\in \ov X^*$ and holomorphic with respect to $\hb$ near $\hb_o$. The action of $\hb\partial_\hb$ is given by $\hb\partial_\hb m=\sum_{p,j}(pf_{p,j}+\hb\partial f_{p,j}\partial\hb) \hb^p \epsilon_{p,j}$.

The condition that $m$ is a local section of $V_{<0}\wt\cM$ is then equivalent to the fact that, for any ${p,j}$, the map $\hb\mto [x\mto xf_{p,j}(x,\hb)\hb^p]$ is a holomorphic function from a neighbourhood of $\hb_o$ to the Banach space of continuous functions on $\ov X$ vanishing at~$0$. It is then clear, from Cauchy's inequality, that $\module{x}\cdot\norme{\hb\partial_\hb m}_{\pi^*h}\leq c\module{x}\cdot\norme{m}_{\pi^*h}$ for some $c>0$.
\end{proof}

It follows that $V_{<0}\wt\cM$ is a locally free $\cO_{\cP^1}$-module with a meromorphic connection~$\nabla$ having regular singularities along its polar locus $\{\hb=0\}\cup (P\times\Omega_0)$. Let us also notice that, as we have seen in \S\ref{subsec:charact}, the monodromy around $\hb=0$ is the identity.

\begin{lemme}\label{lem:RFVneg}
Denote by $V_{<0}\ccM$ the kernel of $\hb\partial_{\hb}$ acting on $V_{<0}\wt\cM_{|\hb=0}[\hbm]$ and by $F_\bbullet V_{<0}\ccM$ the filtration by ``the order of the pole in $\hb$''. Then we have a natural isomorphism $\cO_{\cP^1}\otimes_{\cO_{\PP^1}[\hb]}R_FV_{<0}\ccM\isom V_{<0}\wt\cM$, identifying (the restriction of) $C$ with $R_Fk$.
\end{lemme}

\begin{proof}
Let us fix a local coordinate $x$ on a small disc $X$ centered at a puncture and a local basis of $V_{<0}\wt\cM$ near the point $(0,0)$. The matrix of the connection in this basis can be written as
\[
\wt A(x,\hb)\,\frac{d\hb}{\hb}+\wt B(x,\hb)\,\frac{dx}{\hb x}
\]
with $\wt A(x,\hb),\wt B(x,\hb)$ holomorphic.

The argument used in the proof of Proposition \ref{prop:smoothgraded} extends to the present situation:
\begin{itemize}
\item
For Step \eqref{enum:a}, we can apply exactly the same arguments; at the end, we find that $\wh P$ is a $\hb$-formal horizontal section of a holomorphic bundle with meromorphic connection on $\cX$ having poles along $\{\hb=0\}\cup\{x=0\}$; moreover, the restriction to the curves $x=x_o\neq0$ and $\hb=\hb_o\neq0$ of this bundle with connection has a regular singularity at the origin, as can be seen on its matrix; after Deligne's criterion in \cite{Deligne70}, the meromorphic connection has regular singularities; it follows that any formal solution is convergent. 
\item
Similarly, for Step \eqref{enum:b}, we apply the same argument, where the $\cO_X$-module $\ker\hb\partial_\hb$ is now called $V_{<0}\ccM$ instead of $V$.
\item
The bundle $\cHom_{\cO_\cX}\bigl(\cO_\cX\otimes_{\cO_X[\hb]}R_FV_{<0}\ccM,V_{<0}\wt\cM\bigr)$ is naturally equipped with a meromorphic connection having regular singularities along $\{\hb=0\}\cup\{x=0\}$; we have constructed a germ of horizontal section of this bundle on $X\times\Delta$, where $\Delta$ is a neighbourhood of $0$ in $\Omega_0$, which is an isomorphism of bundles with connection; it extends in a unique way by horizontality as a section on $X^*\times\Omega_0$, and is at most meromorphic along $\{x=0\}\times\Omega_0$; it is in fact holomorphic along $\{x=0\}\times\Omega_0$, being yet holomorphic along $\{x=0\}\times\Delta$; this gives Step \eqref{enum:c}; in particular, we have identified $V_{<0}\ccM=\ker\big[\hb\partial_\hb:V_{<0}\wt\cM_{|\hb=0}[\hbm]\to V_{<0}\wt\cM_{|\hb=0}[\hbm]\big]$ with $V_{<0}\wt\cM/(\hb-1)V_{<0}\wt\cM$;
\item
the equality $C=R_Fk$ holds away from $\{x=0\}$, and both are $L^1_{\loc}$ along $\{x=0\}$ (\cf \cite[\S5.3.c]{Bibi01c} for $C$), thus they coincide as sesquilinear pairing taking values in distributions, hence Step~\eqref{enum:d}.\qedhere
\end{itemize}
\end{proof}

\begin{lemme}
The filtration $F_\bbullet V_{<0}\ccM$ satisfies \eqref{eq:FVneg}.
\end{lemme}

\begin{proof}
Looking back to Step \eqref{enum:b} in the proof of Prop\ptbl \ref{prop:smoothgraded}, the assertion is equivalent to saying that the order of the pole along $\hb=0$ can be computed away from~$\{x=0\}$.
\end{proof}

\subsubsection*{The minimal extension}
We now construct $R_F\ccM$. We continue to work locally near a puncture. We denote $\wt\cM=\cO_\cX[x^{-1}]\otimes_{\cO_\cX}V_{<0}\wt\cM$ and we denote by $\cM$ the $\cR_\cX$-submodule generated by $V_{<0}\wt\cM$ in $\wt\cM$. It is clear that $\wt\cM\subset j_*\cH'$ is stable by $\hb\partial_\hb$. Let us notice moreover that $\cM$ is so, according to the commutation relation $\hb\partial_\hb\partiall_x=\partiall_x(\hb\partial_\hb+1)$. We can therefore define $\wt\ccM$ and $\ccM$ by the procedure above, \ie by taking the kernel of the action of $\hb\partial_\hb$ on the $\hb$-localized modules. We similarly get a filtration $F_\bbullet$ on these modules.

We have then $R_F\wt\ccM=\cO_X[x^{-1},\hb]\otimes_{\cO_X[\hb]}R_FV_{<0}\ccM$ and $R_F\ccM$ is identified with the $R_F\cD_X$-submodule generated by $R_FV_{<0}\ccM$ in $R_F\wt\ccM$. We conclude that $\ccM$ is a minimal extension (\ie $\cD_X$-generated by $V_{<0}\ccM$) and that $F_\bbullet\ccM$ is a good filtration, as $R_F\ccM$ is then clearly $R_F\cD_X$-coherent.

By localization with respect to $x$ of the isomorphism in Lemma \ref{lem:RFVneg}, we get an isomorphism $\cO_\cX\otimes_{\cO_X[\hb]}R_F\wt\ccM\isom\wt\cM$, and by taking the submodules generated by $R_FV_{<0}\ccM$, we get $\cO_\cX\otimes_{\cO_X[\hb]}R_F\ccM\isom\cM$.

\subsubsection*{Extending the sesquilinear pairing}
We have yet obtained a sesquilinear pairing $k$ on $V_{<0}\ccM\otimes_\CC\ov{V_{<0}\ccM}$ with values in $L^1_{\loc}$, such that the sesquilinear pairing $C$ restricted to $R_FV_{<0}\ccM\otimes_{\CC[\hb,\hbm]}\ov{R_FV_{<0}\ccM}$ coincides with $R_Fk$.

Consider local sections $m',m''$ of $F_k\ccM$ and $F_\ell\ccM$. Then $C(m'\hb^k,\ov{m''\hb^\ell})=\hb^k\ov\hb^\ell C(m',\ov{m''})$ belongs to $\Dbh{X}$, and so does $C(m',\ov{m''})$. But $m',m''$ are obtained from $V_{<0}\ccM$ by acting differential operators. It follows that, by sesquilinearity, $C(m',\ov{m''})$ is a distribution on $X$ (\ie does not depend on $\hb$), as this is true on $V_{<0}\ccM$. We denote by $k(m',\ov{m''})$ this distribution. Then, clearly, $k$ is the desired extension, and it satisfies $C=R_Fk$.

\subsubsection*{Description of the filtration $F_\bbullet\ccM$}
In order to end the proof of Proposition \ref{prop:twHodge}, it remains to identify the previously constructed filtration $F_\bbullet\ccM$ with that given by Formula \eqref{eq:FccM}.

Taking the degree $p$ part in $\hb$ of the equality $R_F\ccM=R_F\cD_X\cdot R_FV_{<0}\ccM$ of submodules of $R_F\wt\ccM$ gives
\[
F_p\ccM=\sum_kF_k\cD_X\cdot F_{p-k}V_{<0}\ccM\quad \text{in }\wt\ccM,
\]
which is exactly \eqref{eq:FccM}.\qed

\section{Application to cohomologically tame functions on~smooth~affine~varieties}
\label{sec:Brieskorn}

\subsection{Cohomologically tame functions}
Let $U$ be a smooth complex affine variety of dimension $n+1$ and let $f:U\to\CC$ be a regular function on $U$. We say that $f$ is a \emph{cohomologically tame function} (\cf \cite{Bibi96b}) if there exists a diagram
\[
\xymatrix{
U\ar@{^{ (}->}[r]^-\kappa\ar[dr]_f&X\ar[d]^F\\
&\Afu
}
\]
where $X$ is quasi-projective and $F$ is projective such that, for any $c\in\Afu$, the support of the complex of vanishing cycles $\phi_{F-c}\bR \kappa_*\QQ_U$ is contained in $U$. This implies in particular that the critical points of $f$ in $U$ are isolated. This also implies that the cone $\cC$ of the complex $\kappa_!\QQ_U\to\bR \kappa_*\QQ_U$ is such that $\phi_{F-c}\cC=0$ for any $c\in\CC$, hence each cohomology sheaf of $\bR F_*\cC$ is (locally) constant on $\Afu$. We mainly use this last property, which has been considered by N\ptbl Katz \cite[Th\ptbl14.13.3]{Katz90} in an arithmetic setting and in \cite{Bibi96b} in the complex setting (\cf also \cite{N-Z92,Parusinski95,D-S01} and \cite[\S6.3]{Dimca04}).

\begin{remarque}
There are many examples of cohomologically tame functions. Namely, if $f$ is a polynomial on $\CC^{n+1}$, then it is cohomologically tame with respect to the closure $X$ in $\PP^{n+1}\times\CC$ of the graph of $f$ (contained in $\CC^{n+1}\times\CC$ if and only if it satisfies \emph{Malgrange condition} (\cf \cite{Parusinski95}). Polynomials which are tame in the sense of Broughton \cite{Broughton88} are examples of this kind. Recall also that (Laurent) polynomials which are convenient and nondegenerate with respect to their Newton polyhedra are cohomologically tame.
\end{remarque}

Recall that we denote by $\pQQ_U$ the complex $\QQ_U[\dim U]$. Let us recall some properties of the complex $\bR f_*\pQQ_U$. We denote by $p_1,\dots,p_r$ the critical values of $f$ and by $j:\Afuan\moins\{p_1,\dots,p_r\}\hto\Afuan$ the open inclusion. We will use basic results concerning perverse sheaves and intersection cohomology, for which we refer to \cite{Dimca04} and the references therein. We also use the nearby and vanishing cycle functors $\psi_g,\phi_g$ relative to a function $g$ (see \loccit), and their perverse counterpart $\psip_g=\psi_g[-1]$, $\phip_g=\phi_g[-1]$.

A basic point will be to compare the direct image $\bR f_*\pQQ_U$ 
with the direct image by $F$ of the intersection complex 
$\IC_X(\pQQ)$.

\begin{proposition}\label{prop:imdirminext}
The perverse sheaf $\pcH^0\bR f_*\pQQ_U$ shifted by $-1$ is a sheaf,
with fibre at $c$ equal to $H^n(f^{-1}(c),\QQ)$. Moreover, there is 
an exact sequence
\begin{equation}\tag*{(\protect\ref{prop:imdirminext})$(*)$} \label{eq:shortex}
0\to\cF_1\to\pcH^0\bR F_*\IC_X(\pQQ)\to\pcH^0\bR f_*\pQQ_U\to\cF_2\to0
\end{equation}
in the perverse category, where $\cF_1,\cF_2$ are constant sheaves shifted by $1$. Last, the perverse sheaf $\pcH^0\bR F_*\IC_X(\pQQ)$ decomposes as the direct sum of two perverse sheaves:
\begin{itemize}
\item
$j_*\cV_{!*}[1]$, where $\cV_{!*}$ is the local 
system with fiber the intersection cohomology $\IH^n(F^{-1}(c),\QQ)$ on $\Afu\moins\{p_1,\dots,p_r\}$,
\item
a sheaf supported on $\{p_1,\dots,p_r\}$, each fiber being of finite 
rank.
\end{itemize}
\end{proposition}

\begin{exemple}
Consider the function $f:\CC^{n+1}\to\CC$ given by 
$f(x_0,\dots,x_n)=\sum x_i^2$. We have $r=1$ and $p_1=0$.
\begin{itemize}
\item
If $n=1$, then $\pcH^0\bR f_*\pQQ_U=j_!\pQQ_{\Afu\moins\{p_1,\dots,p_r\}}$, $\pcH^0\bR 
F_*\IC_X(\pQQ)=i_{0,*}\QQ_0$, so 
$\cV_{!*}=0$, $\cF_1=0$ and $\cF_2=\pQQ_{\Afu}$.
\item
If $n=2$, then $\pcH^0\bR f_*\pQQ_U=j_*\cV$, where $\cV$ has rank one, 
$\pcH^0\bR F_*\IC_X(\pQQ)=j_*\cV_{!*}$, where $\cV_{!*}$ has 
rank two, $\cF_1$ has rank one and $\cF_2=0$.
\end{itemize}
\end{exemple}

\begin{proof}
The first point was proved in \cite[Th\ptbl8.1(3)]{Bibi96b}; more precisely, to identify the fibre, we use that, if $i_{c}: F^{-1}(c) \hto X$ denotes the closed inclusion and $\kappa_c:f^{-1}(c)\hto F^{-1}(c)$ denotes the open inclusion, then we have $i_{c}^{-1}\bR\kappa_*\QQ_U=\bR\kappa_{c,*}\QQ_{f^{-1}(c)}$: this is proved in \cite[Prop\ptbl8.3(3)]{Bibi96b}.

For the second point, let us first remark that the intersection complex $\IC_X(\pQQ)$ also satisfies the non-existence of vanishing cycles at infinity. Indeed, cohomological tameness for $f$ is equivalent to the vanishing along $X\moins U$ of the complexes $\phip_{F-c}\pcH^\ell(\bR \kappa_*\pQQ_U)$ for all~$\ell$, if $\pcH^\ell$ denotes the $\ell$-th complex of perverse cohomology. This follows from the commutation of the functor $\phip_{F-c}$ with the functors $\pcH^\ell$ (see \eg \cite[Cor\ptbl10.3.13]{K-S90}, where our notation $\phip$ corresponds to their notation $\phi$).

On the other hand, one has $\IC_X(\pQQ)=\kappa_{!*}\pQQ_U$ and $\kappa_{!*}\pQQ_U$ is a perverse subsheaf of $\pcH^0(\bR \kappa_*\pQQ_U)$. Therefore, $\phip_{F-c}\kappa_{!*}\pQQ_U$ is a perverse subsheaf of $\phip_{F-c}\pcH^0(\bR \kappa_*\pQQ_U)$, hence is zero on $X\moins U$.

In the same way one shows that, for all $c\in\CC$, the cone of the complex $\kappa_{!*}\pQQ_U\to\bR \kappa_*\pQQ_U$ does not have vanishing cycles along $F=c$. This remains true after direct image by the proper morphism $F$, and thus the perverse cohomology sheaves of the cone $\cC$ of the complex $\bR F_*\IC_X(\pQQ)\to\bR f_*\pQQ_U$ are locally constant perverse sheaves (\ie locally constant sheaves shifted by $1$) on $\Afuan$, hence constant.

On the other hand, still considering the vanishing cycles, one shows that the perverse sheaves $\pcH^\ell\bR F_*\IC_X(\pQQ)$ and $\pcH^\ell\bR f_*\pQQ_U=\pcH^\ell\bR F_*\bR\kappa_*\pQQ_U$ are constant if $\ell\neq0$.

In the long exact sequence of perverse cohomology
\begin{multline*}
\cdots\to\pcH^{-1}\bR f_*\pQQ_U\to\pcH^0\cC\\
\to\pcH^0\bR F_*\IC_X(\pQQ)\to\pcH^0\bR f_*\pQQ_U\\
\to\pcH^1\cC\to\pcH^1\bR F_*\IC_X(\pQQ)\to\cdots
\end{multline*}
all the terms but those in the middle are constant sheaves (shifted by 
one), hence we get a short exact sequence \ref{eq:shortex}.

The Decomposition Theorem \cite{B-B-D81} gives the desired 
decomposition. Let us show that the fibre of $\cV_{!*}$ at $c\in\Afu\moins\{p_1,\dots,p_r\}$ is the intersection
cohomology $\IH^n(F^{-1}(c),\QQ)$, and the morphism
$\cV_{!*,c}\to\cV_c$ is the natural restriction morphism
$\IH^n(F^{-1}(c),\QQ)\to H^n(f^{-1}(c),\QQ)$.

To prove this statement, it is enough to prove that, in some neighbourhood of $X\moins U$, we have, for any $c\in\Afuan$, the equality $i_{c}^{-1}\kappa_{!*}\QQ_U=\kappa_{c,!*}\QQ_{f^{-1}(c)}$. Let us use the shifting convention for perverse sheaves. Recall that $\kappa_{c,!*}\pQQ_{f^{-1}(c)}$ is the image (in the perverse category) of the natural morphism $\pcH^0\kappa_{c,!}\pQQ_{f^{-1}(c)}\to\pcH^0\kappa_{c,*}\pQQ_{f^{-1}(c)}$. In some neighbourhood of $X\moins U$, according to the non-existence of vanishing cycles of the complexes involved, we have, setting $\ppi_{c}^{-1}=i_{c}^{-1}[-1]$,
\begin{align*}
\ppi_{c}^{-1}\kappa_{!*}\pQQ_U&=\psip_{f-c}\kappa_{!*}\pQQ_U\quad (\text{no vanishing cycles})\\
&=\image\Big[\psip_{f-c}\pcH^0\kappa_!\pQQ_U\ra\psip_{f-c}\pcH^0\bR\kappa_*\pQQ_U\Big]\\
&=\image\Big[\pcH^0\psip_{f-c}\kappa_!\pQQ_U\ra\pcH^0\psip_{f-c}\bR\kappa_*\pQQ_U\Big]\quad (\psip_{f-c}\text{ is $t$-exact})\\
&=\image\Big[\pcH^0\ppi_{c}^{-1}\kappa_!\pQQ_U\ra\pcH^0\ppi_{c}^{-1}\bR\kappa_*\pQQ_U\Big]\quad (\text{no vanishing cycles})\\
&=\image\Big[\pcH^0\kappa_!\ppi_{c}^{-1}\pQQ_U \ra \pcH^0\bR\kappa_*\ppi_{c}^{-1}\pQQ_U\Big]\quad (\text{\cf \cite[Prop\ptbl8.3(3)]{Bibi96b}})\\
&=\image\Big[\pcH^0\kappa_{c,!}\pQQ_{f^{-1}(c)} \ra \pcH^0\bR\kappa_{c,*}\pQQ_{f^{-1}(c)}\Big]\\
&=\kappa_{c,!*}\pQQ_{f^{-1}(c)}.\qedhere
\end{align*}
\end{proof}

\subsection{Semisimplicity of the Gauss-Manin system}
The Gauss-Manin system $M$ of $f$, defined as $\cH^0f_+\cO_U$, is known to be a regular holonomic $\Clt$-module. Its Laplace transform $\wh M$ is a $\Cltau$-module having thus a regular singularity at $\tau=0$ and a (usually) irregular singularity at $\tau=\infty$, and no other singularity. Then $G\defin\CC[\tau,\tau^{-1}]\otimes_{\CC[\tau]}\wh M$ is a free $\CC[\tau,\tau^{-1}]$-module of rank $\mu$ (sum of the Milnor numbers of the critical points of $f$ on $U$), equipped with a connection $\wh\nabla$ (regular singularity at $\tau=0$, irregular one at $\tau=\infty$, and no other singularity).

\begin{theoreme}\label{th:semisimple}
With these assumptions, the meromorphic bundle with connection $(G,\wh\nabla)$ is semisimple, \ie is a direct sum of simple (that is, irreducible) meromorphic bundles with connection.
\end{theoreme}

R\ptbl Garc\'{\i}a\kern2pt L\'opez informed me that this result has been shown by N\ptbl Katz \cite[Th\ptbl14.13.4(4)]{Katz90} by arithmetical methods.

\begin{proof}
Let $M_{\dag+}$ be the regular holonomic $\CC[t]\langle\partial_t\rangle$-module having $\pcH^0\!\bR F_*\IC_X(\pCC)$ as its de~Rham complex (using the Riemann-Hilbert correspondence on $\PP^1$ and GAGA). From \ref{eq:shortex}, we get an exact sequence
\begin{equation}\label{eq:Mdag}
0\to N_1\to M_{\dag+}\to M\to N_2\to0,
\end{equation}
where $N_1,N_2$ are isomorphic to powers of $\CC[t]$ (with its usual $\Clt$-action). In particular, setting $G_{\dag+}=\CC[\tau,\tau^{-1}]\otimes_{\CC[\tau]}\wh{M_{\dag+}}$, one has
\begin{equation}\label{eq:gdagg}
G_{\dag+}=G.
\end{equation}

Let us notice that if $M_{\dag+}$ is semisimple as a $\Clt$-module, then $\wh{M_{\dag+}}$ is so as a $\Cltau$-module and $\CC[\tau,\tau^{-1}]\otimes_{\CC[\tau]}\wh{M_{\dag+}}$ is so as a meromorphic bundle with connection.

The semisimplicity of $M_{\dag+}$ follows from the Decomposition Theorem for the direct image of an intersection complex by a projective morphism (\cite{B-B-D81} or \cite{MSaito86}, see also \cite{C-M03}) and of the Semisimplicity Theorem \cite{DeligneHII}.
\end{proof}

\subsection{The Brieskorn lattice}
We denote by $G_0$ the Brieskorn lattice associated to $f$ (see \cite{Bibi96b}): by definition, one introduces a new variable $\theta$ and one sets
\[
G_0=\Omega^{n+1}(U)[\theta]\big/(\theta d-df\wedge)\Omega^n(U)[\theta].
\]
It is known that $G_0$ is a free $\CC[\theta]$-module of rank $\mu\defin\dim\Omega^{n+1}(U)/df\wedge\Omega^n(U)$. The multiplication by $f$ on $\Omega^n(U)$ naturally induces a connection on $G_0$ with a pole of order $2$ at $\theta=0$, a regular singularity at $\theta=\infty$, and no other pole.

We have $G=\CC[\theta,\theta^{-1}]\otimes_{\CC[\theta]}G_0$, equipped with the corresponding connection.

\subsubsection*{Comparison with $G_{\dag+,0}$}
After \cite{MSaito86}, the module $M_{\dag+}$ is equipped with a good filtration that we denote by $F_\bbullet^{\rH}M_{\dag+}$, so that $(M_{\dag+},F_\bbullet^{\rH}M_{\dag+})$ underlies a polarizable Hodge Module of weight $\dim U=n+1$ on $\Afu$ (which corresponds, on $\Afu\moins\{p_1,\dots,p_r\}$, to a variation of polarized Hodge structure of weight $n$). We will be interested in the shifted filtration $F_\bbullet M_{\dag+}=F_\bbullet^{\rH}[n+1]M_{\dag+}$, defined by $F_\ell M_{\dag+}=F^{\rH}_{\ell-(n+1)}M_{\dag+}$ for any $\ell\in\ZZ$. We then define $G_{\dag+,0}$ as $G_{\dag+,0}^{(F_\bbullet)}$.

\begin{lemme}\label{lem:gdagg}
We have $G_{\dag+,0}=G_0$.
\end{lemme}

\begin{proof}[Proof of $G_{\dag+,0}\subset G_0$]
By definition of the direct image $f_+$, the module $M$ can be written as $\Omega^{n+1}(U)[\partial_t]\big/(d-df\!\wedge\!\cbbullet\,\partial_t)\Omega^n(U)[\partial_t]$. We denote by $M_0$ the image of $\Omega^{n+1}(U)$ in $M$. Then $G_0$ is the saturation \eqref{eq:satdtm} of $M_0$ by $\ptm$. It is therefore enough to show that the natural map $M_{\dag+}\to M$ of \eqref{eq:Mdag} sends $F_0M_{\dag+}$ into $M_0$ and that $0$ is a generating index for $F_\bullet M_{\dag+}$ (\cf \eqref{eq:defG0F}).

We consider a diagram
\[
\xymatrix{
&\wt X\ar[d]^-{\pi}\ar@/^1.5pc/[dd]^-{\wt F}\\
U\ar@{^{ (}->}[r]^-\kappa\ar@{^{ (}->}[ru]^-{\wt\kappa}\ar[dr]_f&X\ar[d]^(.4)F\\
&\Afu
}
\]
where $\wt X$ is smooth and $\pi$ is projective and birational. The direct image $\cH^0\wt F_+\cO_{\wt X}$ is the $0$-th cohomology of the complex
\[
\bR\Gamma\Big(\wt X,\big(\Omega_{\wt X}^{n+1+\cbbullet}[\partial_t],d-d\wt F\!\wedge\!\cbbullet\,\partial_t\big)\Big)
\]
and we have a natural restriction morphism
\begin{multline*}
\bR\Gamma\Big(\wt X,\big(\Omega_{\wt X}^{n+1+\cbbullet}[\partial_t],d-d\wt F\!\wedge\!\cbbullet\,\partial_t\big)\Big)\to \bR \Gamma\Big(U,\big(\Omega_U^{n+1+\cbbullet}[\partial_t],d-df\!\wedge\!\cbbullet\,\partial_t\big)\Big)\\
=\big(\Omega^{n+1+\cbbullet}(U)[\partial_t],d-df\!\wedge\!\cbbullet\,\partial_t\big)
\end{multline*}
inducing a natural morphism
\begin{equation}\label{eq:wFdagM}
\cH^0\wt F_+\cO_{\wt X}\to M.
\end{equation}
We filter the de~Rham complex by
\begin{equation}\label{eq:deRhamfilt}
F_k\big(\Omega_{\wt X}^{n+1+\cbbullet}[\partial_t],d-d\wt F\!\wedge\!\cbbullet\,\partial_t\big)=\big(F_{k+\cbbullet}\big[\Omega_{\wt X}^{n+1+\cbbullet}[\partial_t]\big],d-d\wt F\!\wedge\!\cbbullet\,\partial_t\big)
\end{equation}
with
\[
F_{k+\ell}\big[\Omega_{\wt X}^{n+1+\ell}[\partial_t]\big]=\sum_{a=0}^{k+\ell} \Omega_{\wt X}^{n+1+\ell}\partial_t^a.
\]
This induces the filtration $F_\bbullet\cH^0\wt F_+\cO_{\wt X}$. We have a similar definition for $U$ instead of $\wt X$, and \eqref{eq:wFdagM} is, by construction, compatible with the filtrations. We remark that $F_0M=M_0$; hence, the morphism \eqref{eq:wFdagM} sends $F_0\cH^0\wt F_+\cO_{\wt X}$ into $M_0$.

By construction, $F_0\cH^0\wt F_+\cO_{\wt X}$ generates $F_\bbullet\cH^0\wt F_+\cO_{\wt X}$ over $\CC[\partial_t]$, as
\[
F_0\big(\Omega_{\wt X}^{n+1+\cbbullet}[\partial_t],d-d\wt F\!\wedge\!\cbbullet\,\partial_t\big)=\Omega_{\wt X}^{n+1}.
\]

In order to conclude, we apply \cite{MSaito86}: setting as above $F^{\rH}_\bbullet=F_\bbullet[-(n+1)]$, $(\cH^0\wt F_+\cO_{\wt X},F^{\rH}_\bbullet)$ underlies a polarizable Hodge module of weight $n+1$, and, by the decomposition theorem of \loccit for $\pi_+\cO_{\wt X}$, $(M_{\dag+},F^{\rH}_\bbullet M_{\dag+})$ is a direct summand of $(\cH^0\wt F_+\cO_{\wt X},F^{\rH}_\bbullet)$. Last, the morphism \eqref{eq:Mdag} is induced by \eqref{eq:wFdagM}.
\end{proof}

\begin{proof}[End of the proof of Lemma \ref{lem:gdagg}]
Once we know that $G_{\dag+,0}$ and $G_0$ are two lattices of $G=G_{\dag+}$ and that the former is included in the latter, the lemma follows from the equality after tensoring with $\CC\lcr\theta\rcr$ over $\CC[\theta]$. From the formal stationary phase formula (see \eg \cite[Prop\ptbl V.3.6]{Bibi00}), this reduces to the formal microlocal equality $(M,F_\bbullet M)^\mu=(M_{\dag+},F_\bbullet M_{\dag+})^\mu$ at each point $p_j$.

In order to obtain such an equality, one interprets both microlocal filtered modules as microlocal direct images by $F$ of $\kappa_+\cO_U$ and $\kappa_{\dag+}\cO_U$ respectively (see \eg \cite[\S11]{Bibi96b}). In such a direct image, no contribution comes from $X\moins U$, and the contributions coming from the critical points of $f$ in $U$ coincide.
\end{proof}

\begin{remarque*}[Birkhoff Problem for the Brieskorn lattice]
It follows from Theorem \ref{th:semisimple} and of a theorem of Bolibrukh and Kostov that Birkhoff's problem for the Brieskorn lattice has a solution. This result can be also obtained in another way, which is more precise: Hodge Theory allows one to apply M\ptbl Saito's criterion \cite{MSaito89} to obtain a remarkable solution to this problem (\cf \cite{Bibi96b,D-S02a}).
\end{remarque*}

\subsection{The sesquilinear pairing}\label{subsec:sesquif}
By Poincar\'e-Verdier duality theorem, we have a natural pairing (that we consider from the sesquilinear point of view)
\[
\rP: \bR f_!\pCC_U\otimes_\CC\ov{\bR f_*\pCC_U}\to\CC_{\Afuan}[2].
\]
As the perverse cohomology sheaves in degree distinct from $0$ of both complexes $\bR f_!\pCC_U$ and $\bR f_*\pCC_U$ are constant shifted by one, their topological Laplace transforms vanish on $\Afuhs$, and the topological Laplace transforms of these complexes take the form $\wh\cV[1]$, if $\wh\cV$ is the sheaf of horizontal sections of $G$ on $\Afuhs$. By topological Fourier transform, we get a sesquilinear pairing \eqref{eq:sheafFLsesqui}:
\[
\wh{\rP}:\wh\cV_{|S^1}[1]\otimes_\CC\ov{\iota^{-1}\wh\cV_{|S^1}[1]}\to\CC_{S^1}[2].
\]

For any integer $\ell$, let us set $\epsilon(\ell)=(-1)^{\ell(\ell-1)/2}$.

\begin{theoreme}\label{th:Brieskorn}
Let us assume that $f$ is cohomologically tame. Then the integrable triple $(\cH,\cH,\wh C)$ associated with $(G_0,G_0,[\sfrac{\epsilon(n+1)}{(2\pi i)^{n+1}}]\cdot\wh{\rP})$ by the twistorization process of Definition \ref{def:twistorization} (\ie $\theta\mto\hb$ and replacement of the standard conjugation with the twistor conjugation), equipped with the sesquilinear duality $\cS=(\id_{\cH},\id_{\cH})$, is an integrable polarized twistor structure of weight~$0$.
\end{theoreme}

\begin{proof}
We will apply Corollary \ref{cor:main} to a suitable variation of polarized Hodge structure.

Let us choose a relatively ample line bundle for $F$. Denote by $PM_{\dag+}$ the corresponding primitive submodule of $M_{\dag+}$ (its fibre at $c\in\Afu\moins\{p_1,\dots,p_r\}$ is the primitive part of $\IH^n(F^{-1}(c),\CC)$). The filtration $F_\bbullet^{\rH}M_{\dag+}$ induces a filtration $F_\bbullet^{\rH}PM_{\dag+}$.

\begin{lemme}\label{lem:PV}
The primitive filtered module $(PM_{\dag+},F_\bbullet^{\rH}PM_{\dag+})$ is a direct summand in $(M_{\dag+},F_\bbullet^{\rH}M_{\dag+})$, the other summand being isomorphic to a constant Hodge module.
\end{lemme}

\begin{proof}
By the decomposition theorem \cite{MSaito86}, we have a Lefschetz decomposition
\[
(M_{\dag+},F_\bbullet^{\rH}M_{\dag+})=(PM_{\dag+},F_\bbullet^{\rH}PM_{\dag+})\oplus L \cH^{-2}F_+\,\kappa_{\dag+}(\cO_U,F_\bbullet^\rH\cO_U)
\]
if $(\cO_U,F_\bbullet^\rH\cO_U)$ is $\cO_U$ with its trivial filtration shifted by $-(n+1)$ and $L$ denotes the Lefschetz operator relative to the chosen relatively ample line bundle. The second summand is isomorphic to $\cH^{-2}F_+\,\kappa_{\dag+}(\cO_U,F_\bbullet^\rH\cO_U)$. As we have seen in the proof of Proposition \ref{prop:imdirminext}, the corresponding perverse sheaf is constant.
\end{proof}

Poincar\'e-Verdier duality also induces a sesquilinear pairing $\rP_{!*}$ on $\pcH^0\bR F_*\kappa_{!*}\pCC_U$, in a way compatible to $\rP$. We also denote by $\rP_{!*}$ the restriction to the primitive part. Then (\cf\cite{MSaito86}) $[\epsilon(n)/(2\pi i)^n]\rP_{!*}$ induces the polarization of the variation of Hodge structure of weight $n$ on $\Afu\moins\{p_1,\dots,p_r\}$ corresponding to $(PM_{\dag+},F_\bbullet^{\rH}PM_{\dag+})$. It also induces the polarization of the Hodge structure corresponding to the punctual components of $(PM_{\dag+},F_\bbullet^{\rH}PM_{\dag+})$. If we replace $F_\bbullet^{\rH}$ with $F_\bbullet^{\rH}[n+1]$, we have to replace $[\epsilon(n)/(2\pi i)^n]\rP_{!*}$ with $(-1)^{n+1}[\epsilon(n)/(2\pi i)^n]\rP_{!*}$ (\cf Remark \ref{rem:twsh}), that is, with $-[\epsilon(n+1)/(2\pi i)^n]\rP_{!*}$.

By Lemma \ref{lem:gdagg}, $G_{\dag+,0}=G_0$, and by Lemma \ref{lem:PV}, $G_{\dag+,0}$ is equal to the lattice associated to $(PM_{\dag+},F_\bbullet PM_{\dag+})$. Moreover, we have $\wh{\rP_{!*}}=\wh{\rP}$. We can apply the Main Theorem, according to Corollary \ref{cor:main} for the non punctual components of $(PM_{\dag+},F_\bbullet PM_{\dag+})$ (for the punctual components, we apply Example \ref{ex:basicsbis}\eqref{ex:basicsbis1}).

Let us now notice that
\[
-\frac{\epsilon(n+1)}{(2\pi i)^n}\cdot\itwopi\wh{\rP}=\frac{\epsilon(n+1)}{(2\pi i)^{n+1}}\,\wh{\rP}.
\]

Theorem \ref{th:Brieskorn} is then a consequence of the Main Theorem and of the previous results.
\end{proof}

\backmatter
\providecommand{\bysame}{\leavevmode ---\ }
\providecommand{\og}{``}
\providecommand{\fg}{''}
\providecommand{\smfandname}{\&}

\end{document}